\documentclass[12pt]{amsart}
\usepackage{amssymb}
\usepackage{amsmath, amsfonts}
\begin{document}


\newtheorem{theorem}{Theorem}[section]
\newcommand{\mar}[1]{{\marginpar{\textsf{#1}}}} 
\numberwithin{equation}{section}
\newtheorem*{theorem*}{Theorem}
\newtheorem{prop}[theorem]{Proposition}
\newtheorem*{prop*}{Proposition}
\newtheorem{lemma}[theorem]{Lemma}
\newtheorem{corollary}[theorem]{Corollary}
\newtheorem*{conj*}{Conjecture}
\newtheorem*{corollary*}{Corollary}
\newtheorem{definition}[theorem]{Definition}
\newtheorem*{definition*}{Definition}
\newtheorem{remarks}[theorem]{Remarks}
\newtheorem*{remarks*}{remarks}

\newtheorem*{not*}{Notation}
\newcommand\pa{\partial}
\newcommand\cohom{\operatorname{H}}
\newcommand\Td{\operatorname{Td}}
\newcommand\Trig{\operatorname{Trig}}
\newcommand\Hom{\operatorname{Hom}}
\newcommand\End{\operatorname{End}}
\newcommand\Ker{\operatorname{Ker}}
\newcommand\Ind{\operatorname{Ind}}
\newcommand\oH{\operatorname{H}}
\newcommand\oK{\operatorname{K}}
\newcommand\codim{\operatorname{codim}}
\newcommand\Exp{\operatorname{Exp}}
\newcommand\CAP{{\mathcal AP}}
\newcommand\T{\mathbb T}

\newcommand\ep{\epsilon}
\newcommand\te{\tilde e}
\newcommand\Dd{{\mathcal D}}

\newcommand\what{\widehat}
\newcommand\wtit{\widetilde}
\newcommand\mfS{{\mathfrak S}}
\newcommand\cA{{\mathcal A}}
\newcommand\maA{{\mathcal A}}
\newcommand\maF{{\mathcal F}}
\newcommand\maN{{\mathcal N}}
\newcommand\cM{{\mathcal M}}
\newcommand\maE{{\mathcal E}}
\newcommand\cF{{\mathcal F}}
\newcommand\maG{{\mathcal G}}
\newcommand\cG{{\mathcal G}}
\newcommand\cH{{\mathcal H}}
\newcommand\maH{{\mathcal H}}
\newcommand\cO{{\mathcal O}}
\newcommand\cR{{\mathcal R}}
\newcommand\cS{{\mathcal S}}
\newcommand\cU{{\mathcal U}}
\newcommand\cV{{\mathcal V}}
\newcommand\cX{{\mathcal X}}
\newcommand\cD{{\mathcal D}}
\newcommand\cnn{{\mathcal N}}
\newcommand\wD{\widetilde{D}}
\newcommand\wL{\widetilde{L}}
\newcommand\wM{\widetilde{M}}
\newcommand\wV{\widetilde{V}}
\newcommand\Ee{{\mathcal E}}
\newcommand{\npartial}{\slash\!\!\!\partial}
\newcommand{\Heis}{\operatorname{Heis}}
\newcommand{\Solv}{\operatorname{Solv}}
\newcommand{\Spin}{\operatorname{Spin}}
\newcommand{\SO}{\operatorname{SO}}
\newcommand{\ind}{\operatorname{ind}}
\newcommand{\Index}{\operatorname{index}}
\newcommand{\ch}{\operatorname{ch}}
\newcommand{\rank}{\operatorname{rank}}
\newcommand{\G}{\Gamma}

\newcommand{\abs}[1]{\lvert#1\rvert}
 \newcommand{\A}{{\mathcal A}}
        \newcommand{\D}{{\mathcal D}}
\newcommand{\HH}{{\mathcal H}}
        \newcommand{\LL}{{\mathcal L}}
        \newcommand{\B}{{\mathcal B}}
        \newcommand{\K}{{\mathcal K}}
\newcommand{\oo}{{\mathcal O}}
         \newcommand{\PP}{{\mathcal P}}
        \newcommand{\s}{\sigma}
        \newcommand{\coker}{{\mbox coker}}
        \newcommand{\p}{\partial}
        \newcommand{\dd}{|\D|}
        \newcommand{\n}{\parallel}
\newcommand{\bma}{\left(\begin{array}{cc}}
\newcommand{\ema}{\end{array}\right)}
\newcommand{\bca}{\left(\begin{array}{c}}
\newcommand{\eca}{\end{array}\right)}
\newcommand{\sr}{\stackrel}
\newcommand{\da}{\downarrow}
\newcommand{\tD}{\tilde{\D}}
        \newcommand{\R}{\mathbf R}
        \newcommand{\C}{\mathbf C}
        \newcommand{\h}{\mathbf H}
\newcommand{\Z}{\mathbf Z}
\newcommand{\N}{\mathbf N}
\newcommand{\tto}{\longrightarrow}
\newcommand{\ben}{\begin{displaymath}}
        \newcommand{\een}{\end{displaymath}}
\newcommand{\be}{\begin{equation}}
\newcommand{\ee}{\end{equation}}
        \newcommand{\bean}{\begin{eqnarray*}}
        \newcommand{\eean}{\end{eqnarray*}}
\newcommand{\nno}{\nonumber\\}
\newcommand{\bea}{\begin{eqnarray}}
        \newcommand{\eea}{\end{eqnarray}}

\newcommand{\e}{\epsilon}

\newcommand{\supp}[1]{\operatorname{#1}}
\newcommand{\norm}[1]{\parallel\, #1\, \parallel}
\newcommand{\ip}[2]{\langle #1,#2\rangle}
\newcommand{\nc}{\newcommand}
\nc{\gf}[2]{\genfrac{}{}{0pt}{}{#1}{#2}}
\nc{\mb}[1]{{\mbox{$ #1 $}}}
\nc{\real}{{\mathbb R}}
\nc{\comp}{{\mathbb C}}
\nc{\ints}{{\mathbb Z}}
\nc{\Ltoo}{\mb{L^2({\mathbf H})}}
\nc{\rtoo}{\mb{{\mathbf R}^2}}
\nc{\slr}{{\mathbf {SL}}(2,\real)}
\nc{\slz}{{\mathbf {SL}}(2,\ints)}
\nc{\su}{{\mathbf {SU}}(1,1)}
\nc{\so}{{\mathbf {SO}}}
\nc{\hyp}{{\mathbb H}}
\nc{\disc}{{\mathbf D}}
\nc{\torus}{{\mathbb T}}
\newcommand{\tk}{\widetilde{K}}
\newcommand{\boe}{{\bf e}}\newcommand{\bt}{{\bf t}}
\newcommand{\vth}{\vartheta}
\newcommand{\CGh}{\widetilde{\CG}}
\newcommand{\db}{\overline{\partial}}
\newcommand{\tE}{\widetilde{E}}
\newcommand{\tr}{\mbox{tr}}
\newcommand{\ta}{\widetilde{\alpha}}
\newcommand{\tb}{\widetilde{\beta}}
\newcommand{\txi}{\widetilde{\xi}}
\newcommand{\hV}{\hat{V}}
\newcommand{\IC}{\mathbf{C}}
\newcommand{\IZ}{\mathbf{Z}}
\newcommand{\IP}{\mathbf{P}}
\newcommand{\IR}{\mathbf{R}}
\newcommand{\IH}{\mathbf{H}}
\newcommand{\IG}{\mathbf{G}}
\newcommand{\IS}{\mathbf{S}}
\newcommand{\CC}{{\mathcal C}}
\newcommand{\CD}{{\mathcal D}}
\newcommand{\CS}{{\mathcal S}}
\newcommand{\CG}{{\mathcal G}}
\newcommand{\CL}{{\mathcal L}}
\newcommand{\CO}{{\mathcal O}}
\nc{\ca}{{\mathcal A}}
\nc{\cag}{{{\mathcal A}^\Gamma}}
\nc{\cg}{{\mathcal G}}
\nc{\chh}{{\mathcal H}}
\nc{\ck}{{\mathcal B}}
\nc{\cl}{{\mathcal L}}
\nc{\cm}{{\mathcal M}}
\nc{\cn}{{\mathcal N}}
\nc{\cs}{{\mathcal S}}
\nc{\cz}{{\mathcal Z}}
\nc{\sind}{\sigma{\rm -ind}}

\newcommand\clFN{{\mathcal F_\tau(\mathcal N)}}       
\newcommand\clKN{{\mathcal K_\tau(\mathcal N)}}       
\newcommand\clQN{{\mathcal Q_\tau(\mathcal N)}}       %
\newcommand\tF{\tilde F}
\newcommand\clA{\mathcal A}
\newcommand\clH{\mathcal H}
\newcommand\clN{\mathcal N}
\newcommand\Del{\Delta}
\newcommand\g{\gamma}

 \title{An analytic approach to spectral flow in von Neumann algebras}

\author{M-T. BENAMEUR}
\address{UMR 7122 du CNRS\\
 Universit\'e de Metz\\
Ile du Saulcy, Metz,
 FRANCE\\
 e-mail: benameur @math.univ-metz.fr\\}
\author{A. L. CAREY}
\address{Mathematical Sciences Institute,
 Australian National University\\
 Canberra ACT, 0200 AUSTRALIA\\
 e-mail: acarey@maths.anu.edu.au\\}
\author{J. PHILLIPS}
\address{Department of Mathematics and Statistics\\
University of Victoria\\Victoria, B.C. V8W 3P4, CANADA\\
e-mail: phillips@math.uvic.ca}
\author{A. RENNIE}
\address{Institute for Mathematical Sciences\\
Department of Mathematics\\
Universitetsparken 5, DK-2100 Copenhagen DENMARK\\
e-mail: rennie@math.ku.dk}
\author{F. A. SUKOCHEV}
\address{School of Informatics and Engineering,
Flinders University\\
Bedford Park S.A 5042 AUSTRALIA\\
e-mail: sukochev@infoeng.flinders.edu.au}
\author{K. P. WOJCIECHOWSKI}
\address{Department of Mathematics, 
IUPUI (Indiana/Purdue),\\
 Indianapolis, IN, 46202-3216, U.S.A. \\
 e-mail: kwojciechowski@math.iupui.edu}
\begin{abstract}
{The analytic approach to spectral flow is about ten years old. 
In that time it 
has evolved to cover an ever wider range of examples. The most critical 
extension was to replace Fredholm operators in the classical sense by 
Breuer-Fredholm operators in a semifinite von Neumann algebra. The latter 
have continuous spectrum  so that the notion of spectral flow turns out to be 
rather more difficult to deal with. However quite remarkably there is a uniform 
approach in which the proofs do not depend on discreteness of the spectrum of 
the operators in question. The first part of this paper gives a brief account of
 this theory extending and refining earlier results. 
It is then applied in the latter parts of the paper to a series of 
examples.
One of the most powerful tools is an integral formula for spectral flow first 
analysed in the classical setting by Getzler and extended to Breuer-Fredholm
operators by some of the current authors.
This integral formula was known for Dirac operators in a variety of forms ever
since the fundamental papers of Atiyah, Patodi and Singer.
One of the purposes of this exposition is to make contact with this early work 
so that one can understand the recent developments in a proper historical 
context. 
In addition we show how to derive these spectral flow formulae
in the setting
of Dirac operators on (non-compact) covering spaces of a compact spin manifold
using the adiabatic method.
This answers a question of Mathai connecting Atiyah's $L^2$-index 
theorem  to our analytic spectral flow.
Finally we relate our work to that of Coburn, Douglas, Schaeffer and Singer
on Toeplitz operators with almost periodic symbol. We generalise their work to 
cover the case of matrix valued almost periodic symbols
on $\R^N$ using some ideas of Shubin. This provides us with an opportunity to 
describe the deepest part of the theory namely the semifinite local index 
theorem in noncommutative geometry. This theorem, which gives a formula for
spectral flow was recently proved by some of the present authors. It provides 
a far-reaching generalisation of the original 1995 result of Connes and 
Moscovici.}
\end{abstract}
        \maketitle
\section{Introduction}

Spectral flow\footnote{MSC Subject classification:
Primary: 19K56, 58J20; secondary: 46L80, 58J30.}
 is normally associated with
paths of operators with discrete spectrum such as Dirac operators on compact
manifolds.\footnote{This research is supported by the ARC 
(Australia), NSERC (Canada) and an early career grant from the University of
Newcastle.}
Even then it is only in the last decade that analytic definitions
have been introduced (previously the definitions were topological).
Recently it has been discovered that if
 one takes an analytic approach to spectral flow then one can
handle examples where the operators may have zero in the continuous spectrum.

The aim of this article is to give a discussion of
spectral flow in as general an analytic setting as is currently feasible.
In fact we consider unbounded operators affiliated
to a semifinite von Neumann algebra and
give examples where the
phenomenon of spectral flow for paths of such operators
occurs quite naturally. 
There has been a lot published recently on this subject,
which is rather technical although the ideas can be explained 
reasonably simply.
This article is thus partly a review of this theory aimed at  exposing
these recent results to a wider audience. 
As the early papers dealt with von Neumann algebras with trivial centre 
(factors) and the more general situation of non-trivial centre was only recently
completely understood we also felt that it was timely to collect the basic
definitions and results in one place. Moreover
we have rounded out the account with some additional new results and
some carefully chosen illustrative examples.

The methods we use are motivated  by noncommutative geometry however
our results may be stated without using that language.
The novel feature of spectral flow for operators affiliated to a
general semifinite von Neumann algebra is that
the operators in question may have zero in their continuous spectrum.
It is thus rather surprising that spectral flow can even be defined in this
situation.

We focus on
spectral flow for a continuous path of self adjoint
unbounded Breuer-Fredholm operators $\{D(s)=D_0+A(s)\}$ for $s\in [0,1]$
in the sense that $A(s)$ is a norm continuous family of
bounded self adjoint operators in a fixed semifinite von Neumann algebra
$\mathcal N$ and $D(s)$ is affiliated to $\mathcal N$ for all $s\in [0,1]$
 (we will elaborate on all of this terminology in subsequent Sections).
We restrict to the paths of bounded
perturbations because the analytic theory is complete and many
interesting examples exist.
The wider question of paths
where the domain and the Hilbert space $H(s)$, on which $D(s)$ is
densely defined,
varies with $s$ is still under investigation
(see the article by Furutani \cite{Fu} for motivation).
 This situation may arise on manifolds with boundary where one varies the
metric and is a difficult problem unless one makes very specific assumptions.
An approach to this question has been introduced by Leichtnam
and Piazza \cite{LP} building on ideas of Dai and Zhang \cite{DZ}
which in turn is based on unpublished work of F. Wu. It works
for Dirac type operators in both the case of closed manifolds and 
the case of (possibly noncompact) covering spaces.
This new notion is that of spectral section.

Spectral sections enable one to define spectral flow
as an invariant in the K-theory of a certain
$C^*$-subalgebra of the von Neumann algebras that we consider in this
article. We have chosen not to discuss it here because, although it
can handle the case where the space $H(s)$ varies with $s$
we feel that the theory is not yet in final form.
Moreover it reduces in the von Neumann setting to Phillips approach.
Another omission is a discussion of the topological meaning of spectral
flow in the general analytic setting.
We refer the reader to the work of Getzler \cite{G}, 
Boo\ss-Bavnbek et al \cite{BLP}, and Lesch \cite{Le}.

While our aim is to put in one place all of the basic ideas we do not 
include complete proofs instead referring where necessary to the literature.
Thus we start with a  summary of Fredholm theory in a general semifinite
von Neumann algebra $\mathcal N$ with a fixed faithful semifinite
trace $\tau$. We refer to such operators as  `$\tau$-Breuer-Fredholm'
because we can trace the origins of the theory to Breuer
\cite{B1,B2} but we need to refine his theory to take account of 
non-uniqueness of the trace $\tau$ on a von Neumann algebra with 
non-trivial centre.
In this setting we discuss
Phillips' analytic approach to spectral flow for paths of bounded
self adjoint Breuer-Fredholm operators in $\mathcal N$.
Then we  include some
simple analytic examples that show the theory is non-trivial.

The theory for paths $\{D(s)\}_{s\in [0,1]}$
of self adjoint unbounded operators proceeds 
via the map $s\to D(s)\to D(s)(1+D(s)^2)^{-1/2}$. When $\{D(s)\}$ is 
a norm-continuous path of perturbations (of the kind considered above)
of $D(0),$ an unbounded self adjoint $\tau$-Breuer-Fredholm operator,
then its image under this map is a continuous path in the space of bounded 
self-adjoint 
$\tau$-Breuer-Fredholm operators \cite{CP1}. Although (in the case 
$\mathcal N=\mathcal B(\HH)$) spectral flow can be defined 
directly for such paths of unbounded operators \cite{BLP}, we can also 
define spectral flow 
in terms of the corresponding path of bounded self-adjoint
operators.

The second half of the paper is about analytic formulae for spectral flow
that have appeared in the literature. After reviewing these formulae
we relate them to classical theory via a study of
spectral flow of  generalised Dirac operators on compact manifolds without 
boundary
and their covering spaces. A question first raised by 
Mathai \cite{M} is settled by
relating spectral flow to the $L^2$ index theorem. 
The deepest result in the theory is the semifinite local index theorem
which we illustrate by application 
to an example of spectral flow for
differential operators with almost periodic coefficients. 
This is inspired by
work of Shubin \cite{Shu} who initiated this line of enquiry.
The  generalisation to semifinite von Neumann algebras
 of the local index theorem of Connes and Moscovici \cite{CM} 
was achieved in papers of some of the present 
authors \cite{CPRS2,CPRS3} and has other interesting applications
(not included here) for example see 
Pask et al \cite{PRe}.

\section{Preliminaries}

\subsection{Notation}

Our basic reference for von Neumann algebras is Dixmier \cite{Dix} where many
 of the concepts we discuss here are described in detail.
For the theory of ideals of compact operators in a semifinite von Neumann
algebra we refer to Fack et al \cite{FK} and Dodds et al \cite{DPSS}.
Throughout this paper we will consider $\mathcal N$, a semifinite von
Neumann algebra (of
type $I_{\infty}$ or $II_{\infty}$ or mixed type) acting on a separable 
Hilbert space
$\cH$.  We will denote by $\tau$ a fixed faithful, normal semifinite
trace on $\mathcal N$ (with the usual normalization if $\cn$ is a 
type $I_{\infty}$ factor).  The
norm-closed 2-sided ideal in $\mathcal N$ generated by the projections
 of finite
trace (usually called $\tau$-finite projections) will be denoted by 
${\mathcal K}_{\tau\mathcal N}$ or just ${\mathcal K}_\mathcal N$ to
lighten the notation.  The quotient
algebra
${\mathcal N}/{\mathcal K}_\mathcal N$ will be
denoted by ${\mathcal Q}_\mathcal N$ and will be called the 
(generalized) Calkin algebra.  We
will let $\pi$ denote the quotient mapping ${\mathcal N}\to 
{\mathcal Q}_{\mathcal N}$.

We will let ${\mathcal F}$ denote the space of all $\tau$-Breuer-Fredholm 
operators in $\mathcal N$, {\it i.e.} ,
$${\mathcal F}=\left\{T\in\mathcal N\mid\pi (T)\hbox{ is invertible in }
{\mathcal Q}_{\mathcal N}\right\}.$$
We denote by ${\mathcal F}^{sa}$ the space of self adjoint operators in 
$\mathcal F.$ The more interesting part of the space of self adjoint 
$\tau$-Breuer-Fredholm operators in $\mathcal N$ will be denoted by 
${\mathcal F}_{*}^{sa}$, {\it i.e.},
$${\mathcal F}_{*}^{sa}=\left\{T\in {\mathcal F}\mid T=T^{*}\hbox{ and }\pi
(T)\hbox{ is neither positive nor negative}\right\}.$$

\subsection{Some history}

For $\mathcal N$ being the algebra of bounded operators on
$\mathcal H$, i.e. the type $
I_{\infty}$ factor case, Atiyah and
Lusztig \cite{APS1,APS3} 
have defined the {\it spectral flow} of a continuous path
in ${\mathcal F}_{*}^{sa}$ to be the number of eigenvalues (counted with
multiplicities) which pass through $0$ in the positive direction minus the
number which pass through $0$ in the negative direction as one
moves from the initial point of the path to the final point.  This
definition is appealing geometrically as an ``intersection number'' and
has been made precise \cite{G,BW,Ph} although it cannot easily 
be generalised beyond the type $I_{\infty}$ factor.
Other motivating remarks may be found
in Boo\ss-Bavnbek et al \cite{BF,BLP}. 
More importantly, there is no obvious
generalization of this definition if the algebra $\mathcal N$
is of type $II_{\infty}$, where
the spectrum of a self-adjoint Breuer-Fredholm operator is not
discrete in a neighbourhood of zero.  J. Kaminker has described this
as the problem of counting ``moving globs of spectrum''.

In his 1993 Ph.D. thesis, V.S. Perera \cite{P1,P2} gave a definition of the
spectral flow of a {\it loop} in ${\mathcal F}_{*}^{sa}$ for a $II_{\infty}$ 
factor, $\mathcal N$.  He showed that
the space, $\Omega ({\mathcal F}_{*}^{sa})$, of loops based at a unitary $
(2P-1)$ in ${\mathcal F}_{*}^{sa}$, is
homotopy equivalent to the space, ${\mathcal F}$, of all Breuer-Fredholm
operators in the $II_{\infty}$ factor, $P\mathcal N P$.  
Since Breuer \cite{B1,B2} 
showed that
the index map ${\mathcal F}\longrightarrow {\bf R}$ classifies the connected 
components of $
{\mathcal F}$, Perera
defines spectral flow as the composition $sf:\Omega ({\mathcal F}_{*}^{
sa})\longrightarrow {\mathcal F}\longrightarrow {\bf R}$ and so
obtains the isomorphism $\pi_1({\mathcal F}_{*}^{sa})\cong {\bf R}$. He also 
showed that this gives
the ``heuristically correct'' answer for a simple family of loops.

While this is an important and elegant result, it has a couple of
weaknesses.  Firstly, since the map $sf$ is not defined directly and
constructively on individual loops it is not clear why spectral flow
is counting ``moving globs of spectrum''.  Secondly, in the
nonfactor setting where the von Neumann algebra may have summands of finite
type the map may not extend to paths which are not loops 
in any sensible way: in a finite algebra (see 5.1)
there can be paths with nonzero spectral flow, but every loop
has zero spectral flow. 

Phillips' approach \cite{Ph,Ph1} is the following.  Let $\chi$ denote the
characteristic function of the interval $[0,\infty )$.  If $\{B_t\}$ is any
continuous path in ${\mathcal F}_{*}^{sa}$, then $\{\chi (B_t)\}$ is a 
$\underline {\hbox{dis}}$continuous path of
projections whose discontinuities arise precisely because of spectral
flow.  For example, if $t_1<t_2$ are neighbouring path parameters and {\bf if}
the projections $P_i=\chi (B_{t_i})$ commute, then the spectral flow from $
t_1$
to $t_2$ should be $\hbox{trace}(P_2-P_1P_2)$ minus $\hbox{trace}(P_1-P_1P_2)$
($=$ amount of nonnegative spectrum
gained minus amount of nonnegative spectrum lost).  However, this is 
clearly the index of the operator
$P_1P_2:P_2(H)\to P_1(H)$.  If these projections do $\underline {\hbox{not}}$ 
commute then one can
still make sense of this index provided $\pi (P_1)=\pi (P_2)$ in the Calkin
algebra.  This notion was called {\it essential codimension} by Brown,
Douglas and Fillmore \cite{BDF} in the type $I_{\infty}$ case and denoted by
$ec(P_1,P_2)$.  Perera \cite{P1,P2} defined the obvious extension of this
concept to $II_{\infty}$ factors and used it to explain why his definition of
spectral flow gives the ``right'' answer in a representative family of
simple loops.  Phillips' \cite{Ph1,Ph}
new ingredient is the fact that the operator
$P_1P_2:P_2(H)\to P_1(H)$ is always a $\tau$-Breuer-Fredholm operator provided
$||\pi (P_1)-\pi (P_2)||<1$.
While Phillips only proved this in the case of a factor, we observed
in Carey et al \cite{CPS2} that
it works for a general semifinite von Neumann algebra.
We will explain the proof in the next Section, and show that the condition
$||\pi(P_1)-\pi(P_2)||<1$ is necessary and sufficient for $P_1P_2$ to be
$\tau$-Breuer-Fredholm.

Since we can (easily) show that the mapping
$t\mapsto\pi\left(\chi(B_t)\right)$ is continuous, we can partition the 
parameter interval
$a=t_0<t_1<\cdots <t_k=b$ so that on each small subinterval the projections
$\pi\left(\chi (B_t)\right)$ are all close.  Letting $P_i=\chi (B_{
t_i})$ for $i=0,1,\cdots ,k$ we then
define:
$$sf\left(\{B_t\}\right)=\sum_{i=1}^k\,\hbox{Ind}(P_{i-1}P_i).$$
With a little effort this works equally well in both the type $I_{
\infty}$ and
$II_{\infty}$ settings and agrees with all previous definitions of spectral
flow where they exist.  A simple lemma is the key to showing that
$sf$ is well-defined and (path-) homotopy invariant.  Defining $\hbox{Hom}
({\mathcal F}_{*}^{sa})$
to be the homotopy groupoid of ${\mathcal F}_{*}^{sa}$,
Phillips proved the following theorem in the case of a factor.
It extends to the general semifinite case \cite{CPRS3}.

\begin{theorem} If $\cn$ is a general
semifinite von Neumann algebra
then $sf$ 
as defined above is a homomorphism from $\hbox{{\rm Hom}}({\mathcal F}_{
*}^{sa})$ to ${\bf R}$ which restricts to an isomorphism of $\pi_
1({\mathcal F}_{*}^{sa})$ with ${\bf Z}$
(respectively ${\bf R}$) when $\mathcal N$ is a factor of type $I_\infty$
 (respectively,
type $II_\infty$).
\end{theorem}

We note that to show that $sf$ is one-to-one on $\pi_1({\mathcal F}_{*}^{
sa})$ one must rely
on Perera's result that $\Omega ({\mathcal F}_{*}^{sa})\simeq {\mathcal F}$.
We also remark that in paragraphs 7, 8 and 9 of the introduction to the
Atiyah-Patodi-Singer paper \cite{APS3}
 the authors appear to be hinting at the existence of a notion
of spectral flow (for paths of self-adjoint Breuer-Fredholm
operators in a $II_{\infty}$ factor)
to be used as a possible tool in an alternate proof of
their index theorem for flat bundles.
In some sense this hope is realised by 
the generalisation \cite{CPRS2,CPRS3} of the Connes-Moscovici local index 
formula to the
semifinite von Neumann algebra setting.

\section{Breuer-Fredholm theory}

The standard references for Breuer-Fredholm operators in a general
semifinite von Neumann algebra are in Breuer \cite{B1,B2}. In earlier work of 
some of the current authors \cite{CPRS3} 
this theory was extended to handle Breuer-Fredholm operators in a 
skew-corner $P\mathcal N Q$ in the general semifinite situation with a fixed 
(scalar) trace $\tau$ in both the bounded and unbounded cases.
All of the expected results hold but their
proofs are a little more subtle. The most difficult case, index theory
for unbounded Breuer-Fredholm operators will not be covered here.
However, in order to handle more cases (including the case of $\tau$-finite
 von Neumann algebras), we allow our operators to vary
within all of $\mathcal F^{sa}$ and not just in $\mathcal F^{sa}_*.$

If $H_1$ is a subspace of $H$, we
denote the projection onto the closure of $H_1$ by $[H_1]$.

\begin{definition} Let $P$ and $Q$ be projections (not necessarily infinite and
not necessarily equivalent) in
 $\cn$ and let $T\in P\cn Q$.  We let $\ker_Q(T)=\ker\left(T_{|Q(H)}\right
)=\ker(T)\cap Q(H).$
The operator $T\in P\cn Q$ is called $(P\cdot Q)$
{\it $\tau$-Fredholm} if and only if

\noindent(1) $[\ker_Q(T])$ and $[\ker_P(T^{*})]$ are $\tau$-finite in $\cn$, 
and\\
(2) there exists a projection $P_1\leq P$ in $\cn$ with $P-P_1$ $\tau$-finite 
in $\cn$ and $P_1(H)\subseteq T(H)$.

\noindent In this case, we define the $(P\cdot Q)${\it -index
of $T$} to be the number:
$$\hbox{Ind}_{(P\cdot Q)}(T)=\tau\left[\ker_Q(T)\right]-
\tau\left[\ker_P(T^{*})\right].$$
\end{definition}

We will henceforth abbreviate this terminology to $\tau$-Fredholm
or sometimes Breuer-Fredholm
and drop the $(P\cdot Q)$
when there is no danger of confusion.
We observe that if $P=Q$ then this is just the definition of
$\tau$-Fredholm used in Phillips et al \cite{PR}
 in the semifinite von Neumann algebra,
$Q\cn Q$, with the trace being the restriction of $\tau$ to $Q\cn Q$.

We summarize the general
situation of $\tau$-Fredholm operators with different domain and 
range \cite{CPRS3}. We re-iterate that the {\bf order} of 
proving the usual results is crucial in developing the skew-corner case, 
as the various 
projections are neither equivalent nor infinite in general.

\begin{lemma}\label{key}
Let $T\in P\mathcal N Q.$ Then,
\noindent
(1) If $T$ is $(P\cdot Q)$-Fredholm, then $T^*$ is $(Q\cdot P)$-Fredholm and
$Ind(T^*)=-Ind(T)$. If $T=V|T|$ is the polar decomposition,
then $V$ is $(P\cdot Q)$-Fredholm with $Ind(V)=Ind(T)$ and $|T|$ is 
$(Q\cdot Q)$-Fredholm of index $0$.\\
(2) The set of all $(P\cdot Q)$-Fredholm operators in $P\mathcal N Q$ is open 
in the norm topology.
\end{lemma}

\begin{definition}
If $T\in P\mathcal N Q$, then a {\bf parametrix} for $T$ is an operator
$S\in Q\mathcal N P$ satisfying $ST=Q+k_1$ and $TS=P+k_2$ where 
$k_1\in\mathcal{K}_{Q\mathcal N Q}$ and $k_2\in\mathcal{K}_{P\mathcal N P}.$
\end{definition}

\begin{lemma}\label{para}
If the usual assumptions on $\mathcal N$ are satisfied, then 
$T\in P\mathcal N Q$ is
$(P\cdot Q)$-Fredholm if and only if $T$ has a parametrix $S\in Q\mathcal N P$.
Moreover, any such parametrix is $(Q\cdot P)$-Fredholm.
\end{lemma}

\begin{prop}\label{productform}
Let $G,P,Q$ be projections in $\mathcal N$ and let $T\in P\mathcal N Q$ be 
$(P\cdot Q)$-Fredholm and 
$S\in G\mathcal N P$ be $(G\cdot P)$-Fredholm, respectively. Then,
$ST$ is $(G\cdot Q)$-Fredholm and $Ind(ST)=Ind(S)+Ind(T).$
\end{prop}

This proof carefully adapts the original ideas of Breuer \cite{B2} 
in a crucial way.
Finally one is easily able to deduce the following expected results.

\begin{corollary}\label{invariance} (Invariance properties of the 
$(P\cdot Q)$-Index)
Let $T\in P\mathcal N Q.$

\noindent(1) 
If $T$ is $(P\cdot Q)$-Fredholm then there exists $\delta>0$ so that
if $S\in P\mathcal N Q$ and $\n T-S\n<\delta$ then $S$ is 
$(P\cdot Q)$-Fredholm and $Ind(S)=Ind(T).$\\
(2) If $T$ is $(P\cdot Q)$-Fredholm and
$k\in P\mathcal K_{\mathcal N}Q$ then $T+k$ is
$(P\cdot Q)$-Fredholm and $Ind(T+k)=Ind(T).$
\end{corollary}

\section{The analytic definition of spectral flow}

\subsection {Essential Codimension}

If $P,Q$ are projections (not necessarily infinite)
in the semifinite von Neumann
 algebra $\cn$ we
wish to define the {\it essential codimension} of $P$ in $Q$ whenever
$||\pi (P)-\pi (Q)||<1$, where $\pi :\cn\to {\mathcal Q}_\cn$
is the Calkin map.  Once we show
that the operator $PQ\in P\cn Q$ is a $\tau$-Fredholm operator in the
sense of Section 3 then we will define the essential codimension
of $P$ in $Q$ to be $\hbox{Ind}(PQ)$. In case $\cn =\mathcal B (\HH)$
a related result to the following lemma appears in Proposition 3.1 of 
Avron et al \cite{ASS} 
where one of their conditions is in terms of essential spectrum.
Our one condition is in terms of the norm,
and the proof is very different.

\begin{lemma}\label{PQ} If $P,Q$ are projections in the semifinite
von Neumann algebra $\cn$ and $\pi :\cn\to
{\mathcal Q}_\cn$ is the Calkin map,
then $PQ\in P\cn Q$ is $(P\cdot Q)-\tau$-Fredholm if and only if
$||\pi (P)-\pi (Q)||<1.$ 
\end{lemma}

\begin{proof} Suppose $||\pi(Q)-\pi(P)||<1$. Then since
$$\hbox{$||\pi (PQP)-\pi (P)||\leq||\pi (Q)-\pi (P)||
<1$}$$
and
$\pi (P)\left(\cn/{\mathcal K}_\cn\right)\,\pi (P)=(P\cn P)/
{\mathcal K}_{P\cn P}$
we see that $PQP$ is a $\tau$-Fredholm operator in $P\cn P$.
Thus,
$\hbox{}\ker_P(QP)\subseteq\ker_P(PQP)$ and so $\left[\ker_P(QP)\right
]\leq\left[\ker_P(PQP)\right]$ where the
latter is a finite projection in $P\cn P$.  Similarly, $\left[\ker_Q(
PQ)\right]$ is a
finite projection in $Q\cn Q$.  Since the range of $PQ$ contains the range
of $PQP$, and since this latter operator is $\tau$-Fredholm in $P\cn
P$,
there is a projection $P_1\leq P$ so that $\tau(P-P_1)<\infty$ and the range of
$P_1$ is contained in the range of $PQ$.  That is, $PQ$ is
$(P\cdot Q)$-Fredholm.

On the other hand, if $PQ$ is $\tau$-Fredholm then $PQP$ is a 
positive $\tau$-Fredholm operator in $P\cn P$. Letting $p=\pi(P)$ and
$q=\pi(Q)$, we see that $pqp$ is an invertible positive operator
in $p{\mathcal Q}_\cn p$ which is $\leq p,$ so $||p-pqp||<1.$
Similarly, $||q-qpq||<1.$ Now,
$$(p-q)^3 =[p-pqp] - [q-qpq]$$
is the difference of two positive operators, so that:
$$-[q-qpq]\leq (p-q)^3\leq [p-pqp].$$
Hence, 
$$||(p-q)^3||\leq Max\{||p-pqp||,||q-qpq||\} < 1.$$
That is,
$$||\pi(P)-\pi(Q)||=||(p-q)^3||^{1/3} <1.$$

\end{proof}

\begin{definition}
If $P$ and $Q$ are projections in  $\cn$ and if
$||\pi (P)-\pi (Q)||<1$ then the {\it essential codimension of $
P$
in $Q$}, denoted $ec(P,Q)$, is the number 
$\hbox{Ind}(PQ)=\hbox{Ind}_{(P\cdot Q)}(PQ)$.  If $P\leq Q$ it is
exactly the codimension of $P$ in $Q$.
\end{definition}
\begin{lemma} If $P_1$, $P_2$, $P_3$ are projections in
$\cn$ and if $||\pi (P_1)-\pi (P_2)||<\frac{1}{2}$
and $||\pi (P_2)-\pi (P_3)||<\frac{1}{2}$
then $ec(P_1,P_3)=ec(P_1,P_2)+ec(P_2,P_3)$.
\end{lemma}
\begin{proof} Since we also have
$||\pi (P_1)-\pi (P_3)||<1$, the terms in the
equation are all defined by Lemma 4.1.  Translating the equation 
into the language
of index and using Lemma 3.1 and Proposition 3.1 we see that it
suffices to prove that
$\hbox{Ind}\left((P_1P_3)^{*}(P_1P_2P_3)\right)=0$.
But,
$$||\pi\left((P_1P_3)^{*}(P_1P_2P_3)\right)-\pi (P_3)||
=||\pi (P_3P_1P_2P_3)-\pi (P_3)||$$
$$ \leq||\pi (P_1P_2)-\pi (P_3)||
\leq||\pi (P_1P_2)-\pi (P_2)|| +||\pi (P_2)-\pi(P_3)||$$
$$\leq||\pi (P_1)-\pi (P_2)|| +||\pi (P_2)-\pi (P_3)|| <1.$$
Thus, there is a compact $k$ in $P_3\cn P_3$ with
$|| P_3P_1P_2P_3+k-P_3|| <1$.  Hence, $\hbox{Ind}(P_3P_1P_2P_
3)=\hbox{Ind}(P_3P_1P_2P_3+k)=0$
as this latter operator is invertible in $P_3\cn P_3$.~~~\end{proof}

\begin{remarks} If $P$ and $Q$ are projections in $\cn$ with
$|| P-Q|| <1$, then $ec(P,Q)=0.$ To see this, note that
$|| PQP-P||\leq|| Q-P|| <1$ so that $PQP$ is invertible in $
P\cn P$ and hence
$\hbox{range}\,\,P\supseteq\hbox{range}\,\,PQ\supseteq\hbox{range}\,\,
PQP=\hbox{range}\,\,P$.  Thus, $\hbox{range}\,\,PQ=\hbox{range}\,\,
P$
and similarly $\hbox{range}\,\,QP=\hbox{range}\,\,Q$ so the $(P\cdot Q)$ 
index of $PQ$ is 0.
\end{remarks}

\subsection{The general definition}

Recall that $\chi =\chi_{[0,\infty )}$ is the characteristic
function of the interval $[0,\infty )$ so that if $T$ is any self-adjoint
operator in a von Neumann algebra ${\mathcal N}$ then $\chi (T)$ is a 
projection in $
{\mathcal N}$.

\begin{lemma} If ${\mathcal N}$ is a von Neumann algebra, ${\mathcal J}$ is a
norm closed 2-sided
ideal in ${\mathcal N}$, $T$ is a self-adjoint operator in ${\mathcal N}$ and $
\pi (T)$ is invertible
in ${\mathcal N}/{\mathcal J}$ (where 
$\pi :{\mathcal N}\to {\mathcal N}/{\mathcal J}$ is the quotient mapping), then
$\chi\left(\pi (T)\right)=\pi\left(\chi (T)\right).$
\end{lemma}

\begin{proof} Since 0 is not in the spectrum of $\pi (T)$, the left hand side 
is a
well-defined element of the $C^{*}$-algebra ${\mathcal N}/{\mathcal J}$.  
Choose $
\epsilon >0$  so that
$[-\epsilon ,\epsilon ]$ is disjoint from $sp\left(\pi (T)\right).$ Let $
f_1\geq f_2$ be the following
piecewise linear continuous functions on {\bf R}:
$$f_1(t)=\left\{\begin{array}{ll} 1 & \mbox{if}\ t\geq0\\ 
\mbox{linear}&\mbox{on}\ [-\epsilon,0]\\ 0 & \mbox{if}\ t\leq-\epsilon
\end{array}\right.,
\qquad f_2(t)=\left\{\begin{array}{ll} 1 & \mbox{if}\ t\geq\epsilon\\ 
\mbox{linear}&\mbox{on}\ [0,\epsilon]\\ 0 & \mbox{if}\ t\leq0
\end{array}\right..$$

Now, $f_1\geq\chi\geq f_2$ on {\bf R}, but all three functions are equal on $
sp\left(\pi (T)\right)$.
Thus,
$$\chi\left(\pi (T)\right)=f_1\left(\pi (T)\right)=\pi\left(f_1(T)\right
)\geq\pi\left(\chi (T)\right)\geq\pi\left(f_2(T)\right)$$
$$=f_2\left(\pi
(T)\right)=\chi\left(\pi (T)\right).$$
Hence,
$\chi\left(\pi (T)\right)=\pi\left(\chi (T)\right).$
\end{proof}

\begin{definition} Let $\cn$ be a semifinite von Neumann algebra
 with fixed semifinite,
faithful, normal trace, $\tau$.  Let ${\mathcal F}^{sa}$ denote the
space
of all
self-adjoint $\tau$-Fredholm operators in $\cn$.  Let $\{B_t\}$  be any
continuous path in ${\mathcal F}^{sa}$ (indexed by some interval $
[a,b]$).  Then $\{\chi (B_t)\}$
is a (generally discontinuous) path of projections in $\cn$.  By
Lemma 4.3 $\pi\left(\chi (B_t)\right)=\chi\left(\pi (B_t)\right)$ and since 
the spectra of $
\pi (B_t)$ are
bounded away from 0, this latter path is continuous.  By compactness
we can choose a partition $a=t_0<t_1<\cdots <t_k=b$ so that for each
$i=1,2,\cdots ,k$
$$||\pi\left(\chi (B_t)\right)-\pi\left(\chi (B_s)\right)||
<\frac{1}{2}\quad\hbox{for all }t,s\hbox{ in }[t_{i-1},
t_i].$$
Letting $P_i=\chi (B_{t_i})$ for $i=0,1,\cdots ,k$ we define the {\it spectral 
flow of
the path $\{B_t\}$} to be the number:
$$sf\left(\{B_t\}\right)=\sum_{i=1}^k\,ec\left(P_{i-1},P_i\right).$$
\end{definition}

To see that this definition is independent of the partition, it suffices
to see that it is invariant under adding a single point to the
partition.  However, this is exactly the content of Lemma 4.2.

\noindent{\bf Remarks} (i) If $\{B_t\}$ is a path in ${\mathcal F}^{sa}$ and 
if $
t\mapsto\chi (B_t)$ is continuous,
then $sf\left(\{B_t\}\right)=0$.
That is, as expected heuristically,
spectral flow can be nontrivial only when the path $t\mapsto\chi (
B_t)$ has
discontinuities.\\
(ii) For $T\in {\mathcal F}^{sa}$, let
$$\hbox{$N(T)=\left\{S\in {\mathcal F}^{sa}\mid||\pi\left(\chi
(S)\right)-\pi\left(\chi (T)\right)|| <\textstyle{\frac{1}{4}}\right
\}$}.$$
Then $N(T)$ is open in ${\mathcal F}^{sa}$ since $S\mapsto\pi\left
(\chi (S)\right)=\chi\left(\pi (S)\right)$ is continuous
on ${\mathcal F}^{sa}$.  Moreover, if $S_1,S_2\in N(T)$,
then by the definition of
spectral flow, all paths from $S_1$ to $S_2$ lying entirely in $N(T)$ have
the same spectral flow, namely, $ec\left(\chi (S_1),\,\chi (S_2)\right).$

\begin{prop} Spectral flow is homotopy invariant, that is, if
$\{B_t\}$ and $\{B_t^{\prime}\}$ are two continuous paths in ${\mathcal F}^{sa}$
 with $B_0=B_0^{\prime}$  and
$B_1=B^{\prime}_1$ which are homotopic in ${\mathcal F}^{sa}$ via a homotopy 
leaving the
endpoints fixed, then $sf\left(\{B_t\}\right)=sf\left(\{B^{\prime}_
t\}\right)$.
\end{prop}

\begin{proof} Let $H:I\times I\to {\mathcal F}^{sa}$ be a homotopy from $
\{B_t\}$ to $\{B^{\prime}_t\}$.  That is, $H$
is continuous, $H(t,0)=B_t$ for all $t$, $H(t,1)=B_t^{\prime}$ for all $
t$,
$H(0,s)=B_0=B^{\prime}_0$ for all $s$, and $H(1,s)=B_1=B_1^{\prime}$ for all $
s$.  By
compactness we can cover the image of $H$ by a finite number of open
sets $\{N_1,\cdots ,N_k\}$ as in Remark~4.2.  The inverse images of these open
sets, $\{H^{-1}(N_1),\cdots ,H^{-1}(N_k)\}$ is a finite cover of $
I\times I$.  Thus, there
exists $\epsilon_0>0$ (the Lebesgue number of the cover) so that any subset
of $I\times I$ of diameter $\leq\epsilon_0$ is contained in some element of 
this finite
cover of $I\times I$.  Thus, if we partition $I\times I$ into a grid of squares 
of
diameter $\leq\epsilon_0,$ then the image of each square will lie entirely 
within
some $N_i$.  Effectively, this breaks $H$ up into a finite sequence of
``short'' homotopies by restricting $H$ to $I\times J_i$ where $J_
i$ are
subintervals of $I$ (of length $\leq\epsilon_0/\sqrt 2$).  These short 
homotopies have
the added property that for fixed $J_i$ we can choose a single partition
of $I$ so that for each subinterval $J_{\ell}$ of the partition, $
H(J_{\ell}\times J_i)$ is
contained in one of $\{N_1,\cdots ,N_k\}$.  By concentrating on 
the {\it i}th `` short
homotopy'' and relabelling $N_1,\cdots ,N_k$ if necessary we can assume $
H$ is such a ``short homotopy.''
By definition, the sum of the spectral flows of the lower paths 
(i.e. along $I\times \{0\}$) is
$sf\left(\{B_t\}\right).$ Since the spectral flows of the vertical paths
(i.e. along $\{t_k\}\times J_i $)
cancel in
pairs, the sum of the spectral flows of the upper paths 
(i.e., along $I \times \{t_1\}$) equals
$sf\left(\{B_t^{\prime}\}\right)$ and hence $sf\left(\{B_t\}\right
)=sf\left(\{B^{\prime}_t\}\right).$
\end{proof}

\noindent{\bf Examples}.
If $\cn$ is a $II_\infty$ von Neumann factor with trace $\tau$ then 
it is well-known (and not
difficult to prove) that $\cn$ contains an abelian von Neumann
subalgebra isomorphic to $L^{\infty}(\hbox{{\bf R}})$ with the property that
the restriction of the trace $\tau$  to $L^{\infty}({\bf R})$ coincides 
with the usual trace on $L^{\infty}({\bf R})$ given
by Lebesgue integration.  We construct our first examples inside this
subalgebra.  Let $B_0$ in $L^{\infty}({\bf R})$ be the continuous function:
$$B_0(t)=\left\{\begin{array}{cc}
1&\hbox{if }t\geq 1,\\
t&\quad \;\hbox{if }t\in [-1,1],\cr
-1& \hbox{if }t\leq -1.\cr
\end{array}
\right.$$
Let $s$ be any fixed real number.  Then for $t\in [0,1]$ let $B_t$ be defined
by $B_t(r)=B_0(r+ts)$ for all $r\in {\bf R}$.  Clearly $\{B_t\}$ is a 
continuous path in
${\mathcal F}_{*}^{sa}$.  Moreover, $\chi (B_t)=\chi_{[-ts,\infty )}$ which 
differs from $
\chi_{[0,\infty )}$ by the
finite projection $\chi_{[-ts,0)}$ if $s>0$  (or, $\chi_{[0,-ts)}$ if $
s<0$).  Thus,
$\pi\left(\chi (B_t)\right)$ is {\it constant} in ${\mathcal
  Q}_{_{}N}$.
 Hence,
$$P_0=\chi (B_0)=\chi_{[0,\infty )},\ P_1=\chi (B_1)=\chi_{[-s,\infty )}
$$
and
$$sf\left(\{B_t\}\right)=ec(P_0,P_1)=\hbox{Ind}(P_0P_1)$$ 
$$=\tau(P_1-P_0P_1)-\tau(P_0-P_0P_1)=s.
$$
We note that for these examples the spectral pictures are constant!
That is, $sp(B_t)=[-1,1]$ for all $t$ and $sp\left(\pi (B_t)\right)
=\{-1,1\}$ for all $t$.  Thus,
one cannot tell from the spectrum alone (even knowing the
multiplicities) what the spectral flow will be.

These examples may seem paradoxical
as there exists a (strong-operator
topology) continuous path of unitaries
$\{U_t\}$ so that $B_t=U_tB_0U_t^{*}$.  However, there cannot exist a
{\it norm}-continuous path of such unitaries as this would imply that the
path $t\mapsto\chi (B_t)$  is a norm-continuous path of projections which it is
not since $||\chi (B_t)-\chi (B_s)|| =1$ if $s\neq t$.

On the other hand, it is not hard to prove that there is a unitary
$U_1$ in $\mathcal N$ so that $B_1=U_1B_0U_1^{*}$.  Since the unitary group
of $\mathcal N$ is connected
in the norm topology we can find in $\mathcal N$ a norm continuous path $\{
U_t\}$
of unitaries for $t\in [1,2]$ so that $U_1$ is as above and $U_2=I$.  Then we
can extend $\{B_t\}$ to a continuous loop based at $B_0$ by defining
$B_t=U_tB_0U_t^{*}$ for $t\in [1,2].$ Since the second half of the loop 
satisfies
$t\mapsto\chi (B_t)$ is norm continuous, its spectral flow is 0 and so
$sf\left(\{B_t\}_{[1,2]}\right)=r$.

When $\mathcal N$ is a type $I_{\infty}$ factor we can use a similar 
construction with
$\ell^{\infty}({\bf Z})$ in place of $L^{\infty}({\bf R})$ to obtain paths 
with any given integer as
their spectral flow.  Of course, these examples will not have a
constant spectral picture.

\begin{remarks} It is clear from the above definition that spectral flow
does not change under reparametrization of intervals and is additive
when we compose paths by concatenation.  Hence, spectral flow
defines a groupoid homomorphism from the homotopy groupoid,
$\hbox{Hom}({\mathcal F}_{*}^{sa})$ to {\bf Z} in the type $I_{\infty}$ factor
case (respectively, to {\bf R} in the type
$II_{\infty}$ factor case).  By the examples just constructed these 
homomorphisms
are surjective in the case of factors, even when restricted to paths based 
at a point $B_0$ in
$F_{*}^{sa}$, {\it i.e.}, $sf:\pi_1({\mathcal F}_{*}^{sa})\to {\bf Z}$ 
(respectively, {\bf R}) is surjective.  To see that
this group homomorphism is one-to-one on a factor requires the
homotopy equivalence ${\mathcal F}_{*}^{sa}\simeq U(\infty )$ in the 
type $I_{\infty}$ factor case or the
homotopy equivalence \cite{P1,P2} 
$\Omega ({\mathcal F}_{*}^{sa})\simeq {\mathcal F}$ in the type 
$II_{\infty}$ factor case:
in fact, both results only need the somewhat weaker result,
$\Omega ({\mathcal F}_{*}^{sa})\simeq {\mathcal F}$.
\end{remarks}

\section{Spectral flow between self adjoint involutions}

We now revisit the special case which is naturally suggested
by the definition of spectral flow.
Choose projections $P,Q\in \mathcal N$ such that $||\pi(P)-\pi(Q)||<1$ so that 
$QP$ is $\tau$-Fredholm. Let $B_0=2Q-1, B_1= 2P-1$ and introduce the path
$B(t)= (1-t)B_0 + tB_1, 0\leq t \leq 1$. One can easily show in this case
that the path $B_t$ consists of Breuer-Fredholm operators. We are interested
in the spectral flow along this path.
By Definition 4.2 it is
equal to the Breuer-Fredholm index of $QP$ in $P\cn Q$. 
By a careful analysis we will explain why this is the right definition.

First notice that $\ker_P(QP)=\ker(Q)\cap \mbox{ran}(P)$ and
$\ker_Q(PQ)= \ker(P) \cap \mbox{ran}(Q)$. A simple calculation also
yields $ \ker_P(QP)\oplus \ker_Q(PQ)\subset \ker(B_0+B_1)$. Conversely
any element $v$ of $\ker(B_0+B_1)$ satisfies
$v=Pv+Qv$ and hence $PQ(Qv)=0$ and  $QP(Pv)v=0$
implying that $ker(B_0+B_1)\subset
\ker_P(QP)\oplus \ker_Q(PQ)$
(note that it is elementary to check that this
is an orthogonal decomposition and in particular
that $\ker_P(QP) \cap \ker_Q(PQ) =\{0\}$).

Consequently to see what happens as we flow along
$B(t), 0\leq t \leq 1$ we initially track what happens
in $$ker(B_0+B_1)=
\ker_P(QP)\oplus \ker_Q(PQ).$$
Now for $v\in \ker_P(QP)$, $B_0v=-v, B_1v=v$ so that
$B(t)v= (2t-1)v$ and spectrum flows from $-1$ to $1$.
Conversely for $v\in \ker_Q(PQ)$ $B_0v=v, B_1v=-v$ and
$B(t)v=(1-2t)v$. Thus we get flow from $1$ to $-1$.

Hence the spectral flow along the path $\{B_t\}$, denoted $sf\{B_t\}$
is the index of $QP:PH\to QH$ as long as we can show that
there cannot be spectral flow coming in some more complex way
from `outside' $\ker(B_0+B_1)$. We analyse this possibility below.

\begin{remarks}
Spectral flow for the path $\{B(t)\}$ actually occurs
at one point, namely $t=1/2$. To see this we note that $B(t)$ has no kernel
for $t\neq 1/2$ and $\ker (B(1/2))= \ker (B_0+B_1)$. The proof of the former
assertion is elementary because if $B(t)v=0$ then $B_0B_1v=-\frac{1-t}{t}v$
so that, taking norms on both sides we deduce that $1-t\leq t$ or $t\geq 1/2$.
Similarly $B_1B_0v=-\frac{t}{1-t}v$ so that again taking norms we obtain
$t\leq 1/2$. Thus there is only a kernel when $t=1/2$.
\end{remarks}

The analysis of this example is helped by the structure of
the algebra generated by $P$ and $Q$. We have:
\begin{lemma}\label{relation}
Let $U$ be the partial isometry in the polar decomposition of $B_0+B_1.$
Then\\
(i) $(B_0+B_1)B(t)=B(1-t)(B_0+B_1)$\\
(ii) $UB(t)=B(1-t)U$ so that  $UB_0=B_1U$
\end{lemma}
\begin{proof}
(i) This is a straightforward calculation.

\noindent(ii) From (i) we get $(B_0+B_1)^2B(t)=B(t)(B_0+B_1)^2$
so that
$$UB(t)|B_0+B_1|=B(1-t)U|B_0+B_1|$$
and hence on $\ker(B_0+B_1)^\perp$ equation (ii) of the lemma holds.
Because $B(t)$ leaves the kernel of $B_0+B_1$ invariant both sides of
(ii) are zero on this kernel proving the result.
\end{proof}

In the type I factor case one can show there is always a gap in the spectrum
of $B_0+B_1$ about zero. This is because on $\ker(B_0+B_1)^\perp$, the operator 
$B_0+B_1$ is boundedly invertible in the type I factor case so there 
can be no spectral
flow on $\ker(B_0+B_1)^\perp.$ We now show that even in a general 
semifinite von Neumann algebra that there can be {\bf no} spectral flow
when $\ker(B_0+B_1)=\{0\}.$ 

\begin{prop} With the above notation, if $\ker(B_0+B_1)=\{0\},$ then
$sf\{B(t)\}=0.$
\end{prop}

\begin{proof} By assumption 
$$\ker(P)\cap \mbox{ran}(Q)=\{0\}=\ker(Q)\cap\mbox{ran}(P).$$
Now $\B=\{1,(Q-P)^2,(Q+P)\}''$ is a commutative von Neumann algebra, 
so that all the spectral projections of $(Q-P)^2$ lie in $\B$. Now 
$\Vert(Q-P)^2\Vert=\Vert Q-P\Vert^2\leq 1$ and by our assumption, $1$ 
is {\bf not} an eigenvalue of $(Q-P)^2$ because
$$(Q-P)^2x=x\Rightarrow P(Q-P)^2x=Px\Rightarrow PQPx=0\Rightarrow(QP)^2x=0$$
$$\Rightarrow(QP)x=0\Rightarrow Px\in \ker(Q)\cap\mbox{ran}(P)=\{0\}.$$
Hence $Px=0$ and similarly $Qx=0$. Thus $x=(Q-P)^2x=0$, and so
 $\chi_{\{1\}}((Q-P)^2)=0$. Now, since $||\pi(P)-\pi(Q)||<1$ the spectral 
projections 
$$ p_n=\chi_{[1-1/n,1]}((Q-P)^2)$$
are $\tau$-finite for large $n$, and in the commutative algebra $\B$. 
Now, the $p_n$ are decreasing to $\chi_{\{1\}}((Q-P)^2)=0$ and so 
$\tau(p_n)\to 0$. Let $\epsilon>0$ and choose $n$ so that 
$\tau(p_n)<\epsilon$, and note that
$$ p_n=\chi_{[-1,-\sqrt{1-1/n}]\cup[\sqrt{1-1/n},1]}((Q-P))$$
so that $p_n$ commute with $Q-P$. Since $p_n\in\B$, it commutes with $Q+P$, 
and so commutes with both $Q$ and $P$!

We now decompose our space with respect to $1=(1-p_n)+p_n$ and note that both 
$p_n(\HH)$ and $(1-p_n)(\HH)$ are left invariant by all the $B_t$. Hence 
the spectral flow will be the sum of the spectral flows on these two subspaces.
Now, since $\tau(p_n)<\epsilon$, the maximum absolute spectral flow on 
$p_n(\HH)$ is $|\tau(p_n)|<\epsilon$.

On the other hand, on $(1-p_n)(\HH)$ we let $Q_n:=(1-p_n)Q(1-p_n)=Q(1-p_n)$ and 
$P_n:=(1-p_n)P(1-p_n)=P(1-p_n)$ so that 
$$B^n_t:=(1-p_n)B_t(1-p_n)=(1-t)Q_n+tP_n\ \ \mbox{on}\ \ (1-p_n)(\HH).$$
Now, $\Vert Q_n-P_n\Vert=\Vert(Q_n-P_n)^2\Vert^{1/2}\leq (1-1/n)^{1/2}$. So 
$B_0^n-B^n_t=2t(Q_n-P_n)$ and so for $t\leq 1/2$
$$\Vert B_0^n-B_t^n\Vert\leq \Vert Q_n-P_n\Vert\leq (1-1/n)^{1/2}
=:1-\delta_n.$$
So when $t\leq 1/2$ we have
$$\s(B_t^n)\subseteq[-1,1]\cap\overline{Ball_{1-\delta_n}(\s(B_0))}
=[-1,-\delta_n]\cup[\delta_n,1].$$
For $t>1/2$ we have $(1-t)\leq 1/2$ and $\Vert B_1^n-B_t^n\Vert\leq 1-\delta_n$
so again $\s(B_t^n)\subseteq[-1,-\delta_n]\cup[\delta,1]$. Hence there can be 
no spectral flow on $(1-p_n)(\HH)$. Finally, since $|sf(B_t)|<\epsilon$ and 
$\epsilon>0$ was arbitrary, $sf(B_t)=0$.
\end{proof}

\subsection{The case of finite von Neumann algebras}

Let us consider the case where the trace $\tau$ on $\mathcal N$ is finite
so that for any two projections $P,Q$ in $\mathcal N$, $P-Q$ is trace class.
 The size of the positive part of the spectrum of $B_0$
is $\tau(Q)$ and the size of the positive part of the spectrum of
$B_1$ is $\tau(P)$. Thus
it is clear that $\tau(P) -\tau(Q)=\tau(P-Q)$ counts the net amount of spectrum
that has moved across zero
as one moves along the path
$B(t), 0\leq t \leq 1$. So the spectral flow is
$$\tau(P-Q)= \frac{1}{2}\tau(B_1-B_0)= 
\frac{1}{2}\tau\left(\int_0^1\frac{d}{dt}B(t)dt\right)
=\frac{1}{2}\int_0^1 \tau(\frac{d}{dt}B(t)) dt.$$
This simple observation should be compared with later formulae
for spectral flow.

Now by Lemma 5.1
there is a partial isometry $U$ with
$UB_0=B_1U$ on $\ker(B_0+B_1)^\perp$ so that
if $R$ is the projection onto this subspace
$\tau [R(B_1-B_0)]=0$. Thus 
as before, to calculate $\tau(B_1-B_0)$,
it suffices to work in $\ker(B_0+B_1)=\ker_Q (PQ)\oplus \ker_P (QP)$ and then
it is clear that on this space
$\tau(P-Q)=\frac{1}{2}\tau(B_1-B_0)$ is the $\tau$ dimension of $\ker_P (QP)$
minus the $\tau$ dimension of $\ker_Q (PQ)$.

\subsection{Example: APS boundary conditions}

Another way of thinking about spectral flow along $\{B(t)\}$
which is familiar
from Atiyah et al \cite{APS3}  
is to relate it to the index of
the differential operator $\frac{\partial}{\partial t} + B(t)$.
We will briefly sketch this connection for our example of involutions.

Let us suppose that there is a path $w(t), 0\leq t\leq 1$
of vectors in $\mathcal H$ such that $w(0)\in \ker(Q)$
and $w(1)\in ran(P)$. That is, $B_0(w(0))=-w(0)$ and $B_1(w(1))=w(1)$ 
so that this path represents some
flow of spectrum across zero along the path $B(t), 0\leq t\leq 1$.
Assume the path is smooth and consider the equation
\be\label{x.x} w'(t) +B(t)w(t)=0.\ee
By restricting our vectors $w$ to lie in $\ker(B_0+B_1)$ we can easily solve 
this equation. First note that from
\ben (B_0+B_1)B(t)=B(1-t)(B_0+B_1)\een
 we see that $B(t)$, for each $t$ leaves
$\ker(B_0+B_1)$ invariant. We know that for $w=w(0)$ in 
$\ker_P(QP)=ker(Q)\cap ran(P)$,
$B(t)w=(2t-1)w$, so (\ref{x.x}) becomes $w^\prime(t)+(2t-1)w(t)=0$ which
has the solution:
$$w(t)= e^{-(t^2-t)} w(0)$$
noting that $w(0)=w=w(1)$ satisfies the boundary conditions.
Similarly if $w\in \ker_Q(PQ)$ one easily constructs
a solution to the adjoint equation
$$-w'(t) +B(t)w(t)=0.$$

Of course these are APS boundary conditions and we are verifying here
that for the differential operator ${\mathcal B}=\frac{\partial}
{\partial t}+B(t)$
with APS boundary conditions the index of $\mathcal B$ is the spectral flow
along the path $B(t), 0\leq t\leq 1$. More precisely
$\mathcal B$ is densely defined on $L^2([0,1],{\mathcal H})$ with domain
the Sobolev space of $\mathcal H$-valued functions on $[0,1]$ with
$L^2$ derivative. Because $B(t)$ leaves $\ker(B_0+B_1)$ invariant
we can solve ${\mathcal B}w=0$ separately on this space and its orthogonal
complement.

Recall Lemma \ref{relation} where $U$, the partial isometry
in the polar decomposition of $B_0+B_1$, gives an isometry from
 $\ker(B_0+B_1)^\perp$ to itself and satisfies $UB(t)=B(1-t)U$. Suppose 
then that
we have ${\mathcal B}w=0$ where $w$ takes its values
in $\ker(B_0+B_1)^\perp$. Then
$v(t)=Uw(1-t)$ satisfies the adjoint equation
${\mathcal B}^*v=0$ with the adjoint boundary conditions
$v(0)\in Q\mathcal H$, $v(1)\in P\mathcal H$. In other words,
each solution of ${\mathcal B}w=0$ has a counterpart solution
$Uw$ of the adjoint equation and vice versa.
Thus, as expected, the net spectral flow on $\ker(B_0+B_1)^\perp$ must be zero.

\section{Spectral flow for unbounded operators}
The framework is that of noncommutative geometry
in the sense of Alain Connes \cite{Co1,Co2,Co3,Co4,Co5}.
However we need to extend this to cover
odd unbounded $\theta$-{\it summable}
or finitely-summable {\it Breuer-Fredholm modules} for a
unital Banach $*$-algebra, $\mathcal A$. These are pairs $({\mathcal N},D)$
where $\mathcal A$ is represented in the semifinite von Neumann algebra 
${\mathcal N}$ with fixed faithful, normal semifinite trace $\tau$ acting on a
Hilbert space, $H$, and $D$ is an unbounded self-adjoint operator on
$H$
affiliated with ${\mathcal N}$ satisfying: $(1+D^2)^{-1}$ is compact
with the additional side condition that either
$e^{-tD^2}$ is trace class for all $t>0$ ($\theta$-summable)
or  $(1+D^2)^{-1/2}\in {\mathcal L}^n$ for all $n>p$
(with $p$ chosen to be the least real number for which this holds)
and
$[D,a]$ is bounded for all $a$ in a dense $*$-subalgebra of $\mathcal A$.
 The condition $(1+D^2)^{-1/2}\in {\mathcal L}^n$
is known as $n$-summability.
An alternative terminology is to refer to $({\mathcal A}, {\mathcal N}, D)$
as a semifinite spectral triple.
The theory of spectral triples in a von Neuman algebra
 was first exposed in Carey et al \cite{CP1} 
and further developed by Benameur et al \cite{BeF}
and some of the present authors \cite{CPS1,CPRS1,CPRS2,CPRS3,Suk}.

If $u$ is a unitary in this dense $*-$subalgebra then
\begin{center}$uDu^{*} = D+u[D,u^{*}] = D+B$\end{center}
where $B$ is a bounded self-adjoint operator in $\mathcal N$.
We say $D$ and $uDu^*$ are {\bf gauge equivalent}.
The path
\begin{center}$D_t^u:=(1-t)\,D+tuDu^{*}=D+tB$\end{center}
is a ``continuous'' path of unbounded self-adjoint Breuer-Fredholm operators.
More precisely,
\begin{center}$F_t^u:=D_t^u\left(1+(D_t^u)^2\right)^{-{1/2}}$\end{center}
is a norm-continuous path of (bounded) self-adjoint $\tau$-Breuer-Fredholm
operators. The {\it spectral flow} along this path $\{F_t^u\}$ (or $\{
D_t^u\}$) is defined using the first Section via
$sf(\{D_t^u\}):=sf(\{F_t^u\})$.
It recovers the pairing of the K-homology class $[D]$ with the K-theory
class [$u$].

We can relate this spectral flow for the path $\{D_t^u\}$ of unbounded
Breuer Fredholm operators to the relative index of two projections
as follows. Let $\tilde F^u_0$ and $\tilde F^u_1$ be the
partial isometries in the polar decomposition of  $F^u_0$ and $F^u_1$
respectively. By convention these extend to unitaries by making them
the identity on $\ker (F^u_0)$ and $\ker (F^u_1)$ respectively.
We introduce the path $\{\tilde F^u_t\}$ where
$\tilde F^u_t=(1-t)\tilde F^u_0 +t\tilde F^u_1$.
We show below that the spectral flow along $\{F_t^u\}$ is in fact equal to
the spectral flow along $\{\tilde F^u_t\}$.





\begin{prop} Let $(\mathcal A,\mathcal N,D)$ be a semifinite spectral triple, 
and let $u\in\mathcal A$ be a unitary.
Then the spectral flow from $D$ to $uDu^*$ is
 $sf(D,uDu^*) = \Ind(PuP),$ where $P:=\chi(D).$
 \label{SFeqIND}
\end{prop}
\begin{proof}
By the definition
$sf(D,uDu^*):=sf(F_D,F_{uDu^*}),$
where $F_D:=D(1+D^2)^{-1/2}.$ With $\tilde F_D$ 
as defined in the paragraph preceding the proposition introduce
the non-negative spectral projection $P$ of $F_D$ 
by $\tilde F_D=2P-1,$ $\tilde F_{uDu^*}=2Q-1=2(uPu^*)-1.$
{\bf If} $\|\pi(P)-\pi(Q)\|< 1$, then by the definition
$$sf(F_D,F_{uDu^*}):=\Ind(PQ)=ec(P,Q).
$$
 To see that $\|\pi(P)-\pi(Q)\|< 1$, we have \cite{CP1}, 
$F_D-F_{uDu^*}\in\mathcal K_{\tau\mathcal N}$ and
 $$\tF_D-F_D=\tF_D(1-|F_D|)=\tF_D(1-|F_D|^2)(1+|F_D|)^{-1}$$
$$=
\tF_D (1+D^2)^{-1} (1+|F_D|)^{-1} \in\mathcal K_{\tau\mathcal N}.$$
Hence,
$$2(P-Q)=\tF_D-\tF_{uDu^*}=(\tF_D-F_D)+(F_D-F_{uDu^*})+(F_{uDu^*}-\tF_{uDu^*})$$
is also $\tau$-compact, and therefore $\|\pi(P)-\pi(Q)\| = 0 < 1.$
By Lemma \ref{PQ} this shows that the operator $PQ$ is
$(P\cdot Q)$-Fredholm.
Hence, using  the formula $\hbox{Ind} (ST)= \hbox{Ind} (S) +
\hbox{Ind} (T)$ from Section 3 above, we obtain
$$sf(D,uDu^*)=\Ind(PQ)=\Ind(PuPu^*)=\Ind(PuP).
$$
\end{proof}

We conclude this Section with a discussion of a theorem of Lesch \cite{L}.
Let $\clA$ be a unital $C^*$-algebra with a faithful tracial state $\tau,$
$(\pi_\tau,\clH_\tau)$ be the GNS representation of $\clA.$
Let $(\clA,\R,\alpha)$ be a $\tau$-invariant $C^*$-dynamical system.
We will identify $\clA$ with its image $\pi_\tau(\clA).$
Let $\clA \times_\alpha \R$ be the crossed product, so it acts on 
$\clH=L^2(\R, \clH_\tau)=L^2(\R)\otimes\clH_\tau.$
So we have representations $\pi$ and $\lambda$ of $\clA$ and $\R$ given as 
follows:
$\pi(a)$ acts on $\xi\in\clH$ by $\pi(a)\xi(s)=\alpha_s^{-1}(a)\xi(s)$ and 
$\lambda_t \xi(s)=\xi(s-t).$
Let $\clN$ be the von Neumann algebra generated by $\clA \times_\alpha \R.$ 
Clearly,
$\lambda=\{\lambda_t\}_{t\in\R}$ is a one-parameter
group of unitaries in $\clN.$ Let $D$ be its infinitesimal generator, 
that is $\lambda_t=e^{-itD},\ t\in\R.$
We have $\pi(\alpha_t(a))=\lambda_t\pi(a)\lambda_{-t},$ which is equivalent 
to $\pi(\delta(a))=2\pi i [D,\pi(a)],$
where $\delta$ is the infinitesimal generator of $\alpha_t$ and 
$a$ is in the domain of $\delta$ which is a dense *-subalgebra
$\clA_0$ of $\clA$.

In this situation a combination of
 Proposition \ref{SFeqIND} and 
earlier work \cite{CPS2} gives the following
index theorem of Lesch \cite{L} and Phillips and Raeburn \cite{PR}.
\begin{theorem} The triple $(\clA_0,\clH,D)$ 
is a $(1,\infty)$-summable semifinite spectral triple
and for any unitary $u\in\clA_0$, the operator $PuP$ is Breuer-Fredholm 
in $P\clN P$, and
$$sf(D,uDu^*)=\Ind(PuP) = \frac 1{2\pi i} \tau(u\delta(u^*)),
$$
where $P=\chi_{[0,\infty)}(D).$
\end{theorem}

In the case when $\clA=C(\bf T),$ $\alpha_t$ is the rotation by an angle $t$ 
and $\tau$ is the arclength integral on $\bf T,$ then modulo some fiddling 
with $\bf R$ vs. $\bf T$,
we infer the classical Gohberg-Krein theorem \cite{GK1}.
\begin{corollary} If $u$ is a unitary in $C(\bf T)$ which is continuously
differentiable then
$$ sf(D,uDu^*)=\Ind(PuP)=
  \frac{-1}{2\pi i}\int_{\bf T}\frac{u^\prime(x)}{u(x)} dx.
$$
\end{corollary}

In the case when $\clA=CAP(\R)=C(\R_B)$ (i.e. $\clA$ is the $C^*$-algebra of 
all uniformly almost periodic functions on $\R,$
which we identify with the $C^*$-algebra of all continuous functions on the 
Bohr compactification $\R_B$),
$\alpha_t$ is a shift by $t$ and $\tau$ is the Haar integral on $\R_B$, 
we immediately infer the Coburn-Douglas-Schaeffer-Singer theorem.
\begin{corollary} If $u$ is a unitary almost periodic continuously 
differentiable 
function then
$$sf(D,uDu^*)=\Ind(PuP)=\lim\limits_{T\to\infty} 
  \frac{-1}{4\pi iT}\int_{-T}^T \frac{u^\prime(x)}{u(x)}dx.
$$
\end{corollary}

In the case when $\clA=C({\bf T}^2),$ $\alpha_t$ is the Kronecker flow given 
by the vector
field $\partial_x + \theta\partial_y $ with an irrational angle $\theta$
and $\tau$ is the Haar integral on ${\bf T}^2$ we get the following 
example \cite{CP1}.
\begin{corollary} For the unitary element $u(z_1,z_2)=z_2$ from $C({\bf T}^2),$ 
we have
$sf(D,uDu^*)=\Ind(PuP)=\theta.
$
\end{corollary}

\section{Fredholm modules and formulae for spectral flow}

In the case of a finite von Neumann algebra with finite trace $\tau$
for any projections $P$ and $Q$, $P-Q$ is trivially trace class.
For a general semi-finite von Neumann algebra $\mathcal N$
arbitrary projections do not satisfy this property.
However summability conditions on $D$ guarantee that
there is a function $f$
such that $f(P-Q)$ is trace class.
In this setting Carey and Phillips \cite{CP2} 
extended results of Avron et al \cite{ASS}. Specifically,
provided $f(1)$ is nonvanishing and $f$ is odd with $f(P-Q)$
trace class, then
$$sf\{B(t)\}=\mbox{Ind} QP= \frac{1}{f(1)}\tau(f(P-Q)).$$
Starting from these results a somewhat lengthy argument produces
general formulae for spectral flow in the case of
$p$-summable and $\Theta$-summable unbounded Fredholm modules \cite{CP1,CP2}
which we will now describe.

\subsection{The Carey-Phillips  Formulae For Spectral Flow}

It was Singer \cite{Singer} 
who suggested, in the 1974 Vancouver ICM, that spectral flow
and eta invariants were given by integrating a one form. The first paper to
systematically exploit this observation was that of Getzler \cite{G}.
Getzler's paper provided the inspiration for the following extensions.
For paths of Dirac type operators formulae analogous to the ones we describe
here go back to the original papers \cite{APS3}, see for example Chapter 8 of 
Melrose \cite{Me}.

Let $(\A,\HH,D)$ be a $\theta$-summable spectral triple.
We focus here on a couple of the main results of Carey et al \cite{CP2}
in the particular case where we compute the spectral flow from $D_0$
to $D_1=uD_0 u^*$. First we have the formula
\ben sf(D,uD u^*)
=\frac{1}{\sqrt{\pi}}\int_0^1\tau(u[D,u^*]e^{-(D+tu[D,u^*])^2})dt.\een
If $(\A,\HH,D)$ is $n$-summable for some $n>1$ then
\ben sf(D,uD u^*)=
\frac{1}{C_{n/2}}\int_0^1\tau(u[D,u^*](1+(D+tu[D,u^*])^2)^{-n/2})dt,\een
with $C_{n/2}=\int_{-\infty}^\infty(1+x^2)^{-n/2}dx$.

In the type I case the theta summable formula  appeared in Geztler \cite{G}.
The proof of this formula in general
uses a result on spectral flow
for {\bf bounded} self adjoint Breuer-Fredholm operators which we will
briefly explain.

\subsection{Paths of unbounded Breuer-Fredholm operators}

Our approach to spectral flow for a path of unbounded self adjoint operators
affiliated to $\mathcal N$ is to introduce the map
$D\mapsto F_D=D(1+D^2)^{-1/2}$.
By Carey et al \cite{CPRS3}, Section 3 if $\{D(t)\} = \{D(0)+A(t)\}$ is a 
path of unbounded self adjoint $\tau$-Breuer-Fredholm operators affiliated
to $\cn$ where $\{A(s)\}$ is a norm continuous path of bounded self adjoint
operators in $\cn,$ then $\{F_{D(t)}\}$ is a continuous path of self adjoint
$\tau$-Breuer-Fredholm operators in $\cn.$ We then define the spectral flow
of the path $\{D(t)\}$ to be the spectral flow of the path $\{F_{D(t)}\},$
and note that in the case $\cn = \mathcal B(\HH)$, by Boo\ss-Bavnbek 
et al \cite{BLP} one can define $sf(\{D(s)\})$ directly.

The principal difficulty introduced by this point of view is that in practice
the map
$D(t)\mapsto F_{D(t)}$ is hard to deal with when it comes to proving
continuity and
differentiability.  One of the main features of earlier work \cite{CP2}
was to surmount this hurdle.

It is easier to deal with the map $s\mapsto (\lambda -D(t))^{-1}$
where $\lambda$ is in the resolvent set of $D(t)$ and to require
continuity of this map into the bounded operators in $\mathcal N$.
This is equivalent to graph norm continuity.
It is  shown that \cite{BLP,Le} for
${\mathcal N}={\mathcal B}({\mathcal H})$
this resolvent map suffices for a definition of spectral flow
for paths of unbounded self adjoint Fredholm operators.
Unfortunately the case of general semifinite
 $\mathcal N$ seems beyond the scope of these methods.

We let ${\mathcal M}_0= \{D=D_0+A\ |\ A\in {\mathcal N}_{sa}\}$
be an affine space modelled on  ${\mathcal N}_{sa}$.
Let $\gamma=\{D_t=D_0+A(t), a\leq t\leq b\}$ be a
piecewise $C^1$ path in  ${\mathcal M}_0$ with
$D_a$ and $D_b$ invertible.
The spectral flow formula
of Getzler \cite{G} when ${\mathcal N}=B(H)$  is
\begin{eqnarray}
\text{sf}(D_a,D_b)&=& -\int_{\gamma}\alpha_\epsilon +
\frac{1}{2}\eta_\epsilon(D_b)-\frac{1}{2}\eta_\epsilon(D_a)\nonumber\\
&=&\sqrt{\frac{\epsilon}{\pi}}\int_a^b\tau\left(\frac{d}{dt}(D_t)
e^{-\epsilon D_t^2}\right)dt+
\frac{1}{2}\eta_\epsilon(D_b)-\frac{1}{2}\eta_\epsilon(D_a).\nonumber
\end{eqnarray}
where $\eta_\epsilon(D)$
are approximate eta invariant correction terms:
$$\eta_\epsilon (D) = \frac{1}{\sqrt{\pi}}\int_{\epsilon}^{\infty}
\tau\left(De^{-tD^2}\right)t^{-1/2}dt$$
and $\alpha_\epsilon$ is a one form defined on ${\mathcal N}_{sa}$,
the tangent space to ${\mathcal M}_0$, via
$\alpha_{\epsilon}(X)=\sqrt{\frac{\epsilon}{\pi}}\tau(Xe^{-\epsilon D^2}).$

 The most general formula in the bounded case deals with a pair
of self-adjoint $\tau$-Breuer-Fredholm operators
$\{F_j,\; j=1,2\}$, joined by
a piecewise $C^1$ path $\{F_t\}$, $t\in [1,2]$
in a certain affine subspace of the
space of all self-adjoint  $\tau$-Breuer-Fredholm operators.
The spectral flow along
such a path
 is given by
$$sf(F_1,F_2)=\frac{1}{C}
\int_1^2\tau\left(\frac{d}{dt}(F_t)|1-F_t^2|^{-r}
e^{-|1-F_t^2|^{-\sigma}}\right)dt
+\gamma(F_2)-\gamma(F_1) $$
where the $\gamma(F_j)$ are eta invariant type correction terms, $C$ is a
normalization constant depending on the parameters $r\geq 0$ and $\sigma\geq1$.
The affine space in which $\{F_t\}$ live is defined in terms of 
perturbations of one fixed $F_0$ and by the
condition that  $|1-F_t^2|^{-r}e^{-|1-F_t^2|^{-\sigma}}$
is trace class \cite{CP2}.
To see one important place where such complicated formulae arise, one takes
the Getzler expression with $\epsilon=1$ and where the endpoints are 
unitarily equivalent so that the end-point
correction terms cancel:
$$\frac{1}{\sqrt{\pi}}\int_a^b\tau\left(\frac{d}{dt}(D_t)
e^{-D_t^2}\right)dt$$
and does the change of variable $F_t=D_t(1+D_t^2)^{-1/2}$. Then,
$(1+D_t^2)^{-1}=(1-F_t^2)=|1-F_t^2|$, and if one is
careless and just differentiates formally (not worrying about the order of the
factors), one obtains the expression:
$$\frac{e}{\sqrt{\pi}}\int_a^b\tau\left(\frac{d}{dt}(F_t)|1-F_t^2|^{-3/2}
e^{|1-F_t^2|^{-1}}\right)dt.$$
While the actual details are much more complicated, this is the heuristic
essence of the reduction of the unbounded case to the bounded case.

The key observation in this approach is a geometric viewpoint due to Getzler.
He noted \cite{G} that
in the unbounded case the integrand in the spectral flow formula is a one 
form on the affine space
${\mathcal M}_0$. This goes back to the observation by Singer that the eta
invariant itself is actually a one form. One may also explain this fact
from our point of view \cite{CP2}. 
In proving the bounded spectral flow formula one uses in
a crucial way the fact that
$D\to \alpha_D$ where $\alpha_D(X)=\tau(Xe^{|1-F_D^2|^{-1}})$
for $X$ in the tangent space to ${\mathcal M}_0$ at $D$ is an exact one form.

\noindent{\bf Example.} We revisit Corollary 6.3.
The straight line path from $D$ to $uDu^*$ is $D_t^u=D+t\theta 1$ 
for $t\in[0,1].$ As $t$ increases from $0$ to $1,$
the spectral subspaces of the operators $D_t^u$ remain the same, but the 
spectral values each increase by $\theta.$ The spectral subspace of $D$
corresponding to the interval $[-\theta,0), \ E=E_{[-\theta,0)},$ is exactly 
the subspace where the spectral values
change from negative to non-negative. By a calculation very similar to the 
example from Section 4.2, the spectral flow of the path $\{D_t^u\}$
is exactly $\tau(E)$ and since $E=\lambda(\hat g)\otimes 1$ 
where $g=\chi_{[-\theta,0)}$ we have
$$  sf(\{D_t^u\}) = \tau(E) = 
  \int\limits_{-\infty}^\infty \chi_{[-\theta,0)}dr = \theta.$$
It is also easy to verify directly that $\theta$ is the Breuer-Fredholm index 
of the operator $T_u:=PuP$ in the II$_\infty$ factor $P\mathcal N P.$
Finally, using the formula of Section 7.1 with $n/2=1$ we calculate:
\bean
  \frac 1\pi \int\limits_{0}^1 \tau\left(\frac d{dt}(D_t^u)
  \left(1+(D_t^u)^2\right)^{-1}\right)dt
  = \frac 1\pi \int\limits_{0}^1 \tau
  \left( \theta\left(1+(D+t\theta)^2\right)^{-1}\right)dt \\
  = \frac \theta\pi \int\limits_{0}^1 
  \left( \int\limits_{-\infty}^\infty \frac 1{1+(r+t\theta)^2} dr\right)dt
  = \frac \theta\pi \int\limits_{0}^1 
  \left( \int\limits_{-\infty}^\infty \frac 1{1+u^2} du\right)dt = \theta
\eean
which gives the expected result $ \Ind(PuP)=sf(\{D_t^u\})=\theta.$

\section{Spectral Flow, Adiabatic Limits and Covering Spaces}

\subsection{Introductory Remarks}

We start with some observations of 
Mathai \cite{M,Ma} on the motivating example
which arises from the fundamental paper of Atiyah \cite{A}.
Assume that $D_0$ is a self adjoint
Dirac type operator on smooth sections of a bundle over an odd dimensional
manifold $M$. We assume that $M$ is not compact but admits a continuous
free action
of a discrete group $G$ such that the quotient of $M$ by $G$ is a compact
manifold. We assume $M$ is equipped with a $G$ invariant
metric in terms of which $D_0$ is defined and from which we obtain an Hermitian
inner product on the sections of the bundle such that $D_0$ is unbounded and
self-adjoint on an appropriate domain.
In this setting the $G$ action on $M$ lifts to an action by unitary operators
on $L^2$ sections of the bundle. The von Neumann
algebra $\mathcal N$ we consider is the commutant of the $G$ action on 
$L^2$ sections \cite{M}.
 Consider a path of the form
$D_t= D_0 +A(t)$
where $A(t)$ is a bounded self adjoint pseudodifferential operator
depending continuously on $t$ in
the norm topology and commuting with the $G$ action.
Then  $D_t$ is self adjoint Breuer-Fredholm operator affiliated to $\mathcal N$, 
and $D_t(1+D_t^2)^{-1/2}$ is a bounded self adjoint Breuer-Fredholm operator
in $\mathcal N$ for each $t\in \R$. 

While $\mathcal N$ is a semifinite von Neumann algebra
it is not in general a factor. There is a natural trace on $\mathcal N$
(considered by Atiyah \cite{A}
in his account of the $L^2$ index theorem) which we now define.
It will be with respect to this
trace that we calculate spectral flow
along $\{D_0+A(t)\}$ (recall that the type II spectral flow depends on the 
choice of trace non-trivially 
when the algebra $\maN$ has non-trival centre). 
On operators with smooth Schwartz
kernel $k(x,y);\ x,y\in M$ the trace $\tau_G$
is given by taking the fibrewise trace of the kernel on the
diagonal $tr (k(x,x))$ and integrating over a fundamental domain for $G$.
This is the natural trace as may be seen
by recognising that the representation of $G$ we obtain
here is quasi-equivalent to the regular representation.
The regular representation
is determined by the standard trace $\tau_0$ which is
given on an element  $\sum \lambda_g g$  of the group algebra by 
$\tau_0(\sum \lambda_g g)
= \lambda_e$ (with the identity being $e\in G$ and the 
$\lambda_g\in \C$).

The analysis of spectral flow traditionally proceeds by replacing $M$ by 
$M\times S^1$
or $M\times [0,1]$ and considering the Dirac operator on this even dimensional
manifold as in Mathai \cite{M}.  
On covering spaces it
is believed by the experts that one should be able to recover 
an analytic spectral flow formula however
one cannot easily extract a proof from the literature. 
The argument we present in this Section shows how a special case of
the  spectral flow formulae 
discussed in the previous Section
arises naturally from adiabatic limit ideas due originally to 
Cheeger \cite{Ch87}.

\subsection{Easy Adiabatic Formula in even dimensions (EAF)}

We use
an adiabatic process, which leads to the formula for the leading 
term in the expansion of the (difference of) heat kernels. 
This is part of an IUPUI preprint \cite{KPW1} which was never submitted 
for publication.
We discuss the simplest possible case (of a compatible Dirac 
operator with coefficients in an auxiliary bundle) in order 
to avoid use of elliptic estimates as was done in 
Cheeger \cite{Ch87}. 
Therefore our main tool is Duhamel's principle (the expansion 
of the heat kernels of the Dirac Laplacians with respect to 
the perturbation terms \cite{BW,McKS}) 
and we call our result the {\it EAF = Easy Adiabatic Formula}. 
(Let us point out that it is not difficult to imitate Cheeger's proof 
and obtain the {\it EAF} in complete generality 
i.e. for the family of Dirac operators with varying first order part.) 
We present a proof in the case of a closed manifold $M$ and later 
on outline why our argument holds in the case of a continuous free action 
of a discrete group. 

Let $B : C^{\infty}(M;S) \to C^{\infty}(M;S)$ denote 
a compatible Dirac operator acting on sections of a bundle 
of Clifford modules $S$ over a closed, odd-dimensional 
manifold $M$ \cite{BW}. Introduce an auxiliary hermitian vector bundle 
$E$ with hermitian 
connection $\nabla$, and the operator
$B_0 = B \otimes_{\nabla}Id_E$
(see Palais \cite{Pa}, Chapter IV). Let $g : E \to E$ denote 
a unitary bundle automorphism, then we can introduce 
the operator
$$B_1 = gB_0g^{-1} = (Id \otimes g)
(B \otimes_{\nabla}Id_E)(Id \otimes g^{-1}) \ \ .$$
The difference $T = B_1 - B_0 = [g,B]g^{-1}$
is a bundle endomorphism and we want to present a formula 
for the spectral flow of the family 
\begin{equation}\label{e:g1}
\{B_u = B_0 + uT\}_{0 \le u \le 1}
\ \, .
\end{equation}

Spectral flow is a homotopy 
invariant, so we can restrict ourselves to the study 
of the spectral flow of a smooth family 
of self-adjoint operators over $S^1$. 
We introduce a smooth cut-off function 
$\alpha : \R \to \R$, such that
$$\alpha(u) \ = 
 \ \begin{cases} & 0 
\quad \ \text{if}
 \ \ \ u \le 1/4 \ \ , \\
    & \ 1  \quad \ \text{if}
 \ \ 3/4 \le u \ \ . \end{cases}$$
We may also assume that there exists a positive 
constant $c$ , such that
\begin{equation}\label{e:1}
\biggm|\frac{d^k\alpha}{du^k}\biggm| \le c{\cdot}u
\ \,
\end{equation}
for $0\le u \le 1$ , $k = 0,1,2$ .  Now, 
we consider the family
\begin{equation}\label{e:2}
\{B_u = B_0 + \alpha(u)T\}
\ \, ,
\end{equation}
which in an obvious way provides us with 
a family of operators on $S^1$. 
We may also consider the corresponding 
operator $\Dd = \partial_u + B_u$ on the closed manifold 
$N = S^1 \times M$
where $B_u$ is given by the formula (\ref{e:2}). 
The operator $\Dd$ acts 
on sections of the bundle 
$[0,1] \times S\otimes E/_{\cong}$ , where 
the identification is given by
$$(1,y;g(y)w) \cong (0,y;w) \ \ \text{where} \ \ 
w \in S_y \otimes E_y \ \ .$$

\begin{theorem}\label{t:1}
The following formula holds for any $t > 0$
\begin{equation}\label{e:3}
index \ \Dd = sf\{B_u\} = \sqrt{\frac{t}{\pi}}{\cdot}\int_0^1
Tr_M{\dot B}_ue^{-tB_u^2}du
\ \, ,
\end{equation}
where as usual ${\dot B}_u = \frac{dB}{du}$ .
\end{theorem}

The first equality in (\ref{e:3}) goes back 
to the original Atiyah-Patodi-Singer 
paper \cite{APS3}. They also 
proved a formula
$$sf\{B_u\} = \int_0^1{\dot \eta}_udu \ \ ,$$
where $\eta_u$ denotes the $\eta$-invariant of 
the operator $B_u$ . The equality 
$$sf\{B_u\} = \sqrt{\frac{t}{\pi}}{\cdot}
\int_0^1Tr_M{\dot B}_ue^{-tB_u^2}du$$
and more (see below) was proved by the last named author 
around 1991 and published 
in one of his IUPUI preprints \cite{KPW1}. 
The formal paper, which was supposed to contain 
a new discussion of {\it Witten's Holonomy Theorem} 
never appeared and the result  
eventually resurfaced in the paper 
by Getzler \cite{G}. We have to mention that it is not difficult
to manufacture a more straightforward argument to prove the second
equality in (8.4). Here we prove a stronger result that provides the
``adiabatic'' equality on the level of heat kernels, from which (8.4)
and other results not covered in the current paper follow.
 The original proof was a
``toy"  model for a simplified version of Cheeger's 
proof of {\it Witten's Holonomy Theorem} \cite{Ch87} 
and proves the corresponding adiabatic equality on the level 
of the kernels of the corresponding heat operators. 
Therefore we call this equality the 
{\it EAF} = {\it Easy Adiabatic Formula}.

The {\it EAF} is obtained by applying the adiabatic 
process on $N$ in the normal direction 
to a fibre $M$ . We replace the product Riemannian 
metric $g = du^2 + g_M$ by a new metric
\begin{equation}\label{e:4}
g_{\e} = \frac{du^2}{\e^2} + g_M = dv^2 + g_M
\ \, ,
\end{equation}
and let the positive parameter $\e$ run to $0$ . 
The corresponding operator $\Dd_{\e}$ has the following 
representation
\begin{equation}\label{e:5a}
\Dd_{\e}(v,y) = \partial_v + B_{\e v}(y)
\ \, .
\end{equation}
In the equality (\ref{e:5a}) we use a new normal coordinate 
$v = \frac{u}{\e}$ ($y \in M$). The operator $\Dd_{\e}$ 
lives on the manifold $N_{\e} = S^1_{\e} \times M$, 
where $S^1_{\e}$ denotes the circle of length $\frac {1}{\e}$. 
Both $index$ and $sf$ do not change their values under 
the deformation, hence we have the equality 
$$index \ \Dd = index \ \Dd_{\e} = 
sf\{B_{\e v}\} = sf\{B_u\} \ \ .$$

The famous McKean-Singer equality expresses the index 
in terms of the kernels of the heat operators
$$index \ \Dd = index \ \Dd_{\e} = 
Tr \ e^{-t\Dd_{\e}^*\Dd_{\e}} - 
Tr \  e^{-t\Dd_{\e}\Dd_{\e}^*} \ \ ,$$
for fixed $t>0$ as $\e \to 0$. 
Let $k_{\e}(t;(v_1,y_1),(v_2,y_2))$ ($v_i$ 
is the coordinate on $S^1_{\e}$  and $y_i \in M$), 
denote the kernel of the operator 
$e^{-t\Dd^*\Dd_{\e}} - e^{-t\Dd_{\e}\Dd_{\e}^*}$. 
This is the difference of the heat kernels, which are pointwise 
bounded for fixed time $t$ (see Proposition \ref{p:1}), which 
expands into expansion with respect to the parameter $\e$. 
The contributions to the leading term from the kernel 
of the operator $\Dd^*\Dd$ and the kernel of the operator 
$\Dd\Dd^*$ cancel each other, hence this term is equal to $0$. 
Theorem \ref{t:2} provides the formula for the second term 
in the expansion. It follows that the kernel 
$k_{\e}(t;(v_1,y_1),(v_2,y_2))$ at the given point is of 
size of $\e$  and the volume of the manifold $N_{\e}$ is equal 
to $\frac {1}{\e} {\cdot}vol(M)$ hence at the end we obtain 
a finite limit
$$\lim_{\e \to 0} \int_{N_{\e}} 
tr \ k_{\e}(t;(v,y),(v,y))dvdy$$ 
equal to (\ref{e:3}).

Now, we present the kernel on $N_{\e}$, 
which gives the leading term in the expansion 
of $k_{\e}(t;(v_1,y_1),(v_2,y_2))$. 
Let us fix $v_0$ the value of the normal 
coordinate and let $e_{v_0}(t;(v_1,y_1),(v_2,y_2))$ 
denote the kernel of the heat 
operator $e^{-t B^2_{\e v_0}}$. We also introduce 
$\e\alpha'(\e v_0)Te_{v_0}$ , kernel of the operator 
$\e\alpha'(\e v_0)Te^{-t B^2_{\e v_0}}$. 
To get the final product we have to take the convolution 
of kernels. If $k_1$, $k_2$ denote two time-dependent 
operators with smooth kernels on $M$, then 
$k_1*k_2(t) = \int_0^tk_1(s)k_2(t-s)ds$ and 
on the level of the kernels we have the equality
$$k_1*k_2(t;y_1,y_2) = \int_0^tds
\int_Mdz \ k_1(s;y_1,z)k_2(t-s;z,y_2) \ \ .$$
We introduce the kernel
\begin{equation}\label{e:5}
\Ee_{v_0}(t;(v_1,y_1),(v_2,y_2)) = 
2\e\alpha'(\e v_0)e_{\partial_v}(t;v_1,v_2)
(e_{v_0}*Te_{v_0})(t;y_1,y_2),
\end{equation}
where $e_{\partial_v}(t;v_1,v_2)$ denotes the kernel of 
the $1$-dimensional heat operator defined on $\R$ by 
the operator $-\partial_v^2$. Let us also point out 
that $\e\alpha'(\e v_0)T$ is simply equal 
to $\e {\dot B}$ (at $v_0$), where $dot$ denotes 
the derivative with respect to $u$-variable, 
so it is the operator 
$\frac{d}{dv} B_{\e v}\biggm|_{v=v_0}.$ 

At last, we are ready to formulate the {\it EAF}
\begin{theorem}\label{t:2}
For any $t>0$ there exists $\e_0$ and a constant $c > 0$ 
such that for any $0 < \e < \e_0$
\begin{equation}\label{e:6}
\frac {1}{\e}{\cdot}\|k_{\e}(t;(v_0,y_1),(v_0,y_2)) 
- \Ee_{v_0}(t;(v_0,y_1),(v_0,y_2))\| 
\le c{\cdot}\sqrt{\frac {\e}t}
\ \, .
\end{equation}
\end{theorem}

\begin{remarks}\label{r:1}
(1) It has been already pointed out that we only 
present the proof of {\it EAF} for the family (\ref{e:2}). 
The method we use works for any family $\{B_u\}_{0\le u \le 1}$, 
such that for any $0 \le u \le 1$ the difference $B_u - B_0$ is 
an operator of order $0$ and
$$B_1 = gB_0g^{-1} \ \ 
{\text and} \ \ \|B_u - B_0\| \le c{\cdot}u \ \ .$$
If the difference betweeen the operators $B_u$ is 
a $1st$ order operator then we have to follow a more 
complicated version of the argument as presented in 
Cheeger's work \cite{Ch87}. 

(2) The proof we use allows us to replace $\sqrt{\e}$, 
which appears on the right side of (\ref{e:6}), by 
$\e^r$ , for any $0< r < 1$. 
\end{remarks}

\begin{proof}
{\it (2nd Part).} The most technical part of the proof 
is presented in the next subsection. There we use 
the Duhamel's Principle to obtain the equality
$$e^{-t\Del_{1,{\e}}} - e^{-t\Del_{2,{\e}}} = 
\int_0^t \ e^{-s\Del_{1,{\e}}}
(\Del_{2,{\e}} - \Del_{1,{\e}})e^{-(t-s)\Del_{2,{\e}}}ds$$
$$=\int_0^t[e^{-s\Del_0}(\Del_{2,{\e}} - \Del_{1,{\e}})e^{-(t-s)\Del_0}
+(e^{-s\Del_{1,\e}} - e^{-s\Del_0})
(\Del_{2,{\e}} - \Del_{1,{\e}})e^{-(t-s)\Del_{2,{\e}}}]ds $$
$$+ \int_0^te^{-s\Del_0}(\Del_{2,{\e}} - \Del_{1,{\e}})
(e^{-(t-s)\Del_{2,{\e}}} - e^{-(t-s)\Del_0})ds \ \ ,$$
and we show that the kernels of the second 
and third terms on the right side are point-wise 
at most of the size $\e^{\frac 32}$. 
Hence, we only have to study the first term and show that 
it gives the kernel (\ref{e:5}) as $\e \to 0$. 

The operator
\begin{equation}\label{e:abc}
\int_0^te^{-s\Del_0}(\Del_{2,{\e}} - 
\Del_{1,{\e}})e^{-(t-s)\Del_0}ds
\ \, .
\end{equation}
can be represented as
$$\int_0^te^{-s(-\partial_v^2)}
e^{-sB_{\e v_0}^2}(2\e\alpha'(\e v)T)
e^{-(t-s)(-\partial_v^2)}
e^{-(t-s)B_{\e v_0}^2}ds =$$
$$\int_0^te^{-s(-\partial_v^2)}(2\e\alpha'(\e v))
e^{-(t-s)(-\partial_v^2)}
e^{-sB_{\e v_0}^2}Te^{-(t-s)B_{\e v_0}^2}ds \ \ .$$
First, we study $l_{\e}(s,t;v_1,v_2)$ the kernel of 
the part which acts in the normal direction 
$$e^{-s(-\partial_v^2)}(2\e\alpha'(\e v))
e^{-(t-s)(-\partial_v^2)} \ \ .$$
We have
$$l_{\e}(s,t;v_0,v_0) = \int_{-\infty}^{+\infty}
\frac {e^{-\frac {(v_0-z)^2}{4s}}}{\sqrt{4\pi s}}
(2\e\alpha'(\e z))
\frac {e^{-\frac {(v_0-z)^2}{4(t-s)}}}{\sqrt{4\pi (t-s)}}dz \ \ .$$
To simplify, we assume $v_0 = 0$ (this will be justified 
in the next subsection). We obtain
$$l_{\e}(s,t;v_0,v_0) = \frac 1{2\pi}\int_{-\infty}^{+\infty}
\frac {e^{-\frac {tz^2}{4s(t-s)}}}{\sqrt{s(t-s)}}
(\e\alpha'(\e z))dz \ \ .$$
Now, we only have to show that
$$\frac 1{2\pi}\int_{-\infty}^{+\infty}
\frac {e^{-\frac {tz^2}{4s(t-s)}}}{\sqrt{s(t-s)}}
(\e\alpha'(\e z) - \e\alpha'(0))dz| \le 
c\frac {\e^{\frac 32}}{\sqrt{t}} \ \ .$$ 
This is the place where we use (\ref{e:1}). 
We apply here the special case of a trick used 
in the next subsection. 

First, an elementary 
computation shows that
$$\int_{|z|> \frac 1{\sqrt{\e}}}
\frac {e^{-\frac {tz^2}{4s(t-s)}}}{\sqrt{s(t-s)}}
(\e\alpha'(\e z) - \e\alpha'(0))dz|$$
is exponentially small with respect to $\e$ . We have   
$$|\int_{|z|> \frac 1{\sqrt{\e}}}
\frac {e^{-\frac {tz^2}{4s(t-s)}}}{\sqrt{s(t-s)}}
(\e\alpha'(\e z) - \e\alpha'(0))dz| \le 
c\frac {\e^{\frac 32}}{\sqrt{t}}| \le 
\e c_1{\cdot}\int_{|z|> \frac 1{\sqrt{\e}}}
\frac {e^{-\frac {tz^2}{4s(t-s)}}}{\sqrt{s(t-s)}}dz$$
$$=\frac 2{\sqrt{t}}{\cdot}
\int_{r> \sqrt{\frac t{4s(t-s)\e}}}e^{-r^2} \le
\frac 2{\sqrt{t}}{\cdot}e^{-\frac 1{4t\e}} \le 
c_1e^{-\frac {c_2}{t\e}} \ \ .$$

It follows now from (\ref{e:1}) that
$$|\alpha'(\e z) - \alpha'(0)| \le c_3 \sqrt{\e}$$
for $z \le \frac 1{\sqrt{\e}}$. This gives
$$|\int_{|z|< \frac 1{\sqrt{\e}}}
\frac {e^{-\frac {tz^2}{4s(t-s)}}}{\sqrt{s(t-s)}}
(\e\alpha'(\e z) - \e\alpha'(0))dz| \le 
c_3\e^{\frac 32}{\cdot}\int_{-\infty}^{+\infty}
\frac {e^{-\frac {tz^2}{4s(t-s)}}}{\sqrt{s(t-s)}} \le 
c_4\frac {\e^{\frac 32}}{\sqrt{t}} \ .$$
We see that up to a term of order 
$\frac {\e^{\frac 32}}{\sqrt{t}}$ the kernel 
$l_{\e}(s,t;v_1,v_2)$ (for $v_1=v_2=v_0$) is equal to
$$\e\alpha'(v_0)\frac 1{2\pi}\int_{-\infty}^{+\infty}
\frac {e^{-\frac {tz^2}{4s(t-s)}}}{\sqrt{s(t-s)}} = 
\frac {\e\alpha'(v_0)}{\sqrt{t}} \ \ ,$$
and  $l_{\e}(s,t;v_1,v_2)$ can be replaced by 
the kernel of the operator
$$2\e\alpha'(v_0)e^{-t(-\partial_v^2)} \ \ .$$

As a result of these estimates the operator (\ref{e:abc}) 
can be replaced by
$$2\e\alpha'(\e v_0)e^{-t(-\partial_v^2)}\int_0^t
e^{-sB_{\e v_0}^2}Te^{-(t-s)B_{\e v_0}^2}ds \ \ ,$$
which has kernel equal to
$$2\e\alpha'(\e v_0)e_{\partial_v}(t;v_1,v_2)\int_0^tds
\int_Me_{v_0}(s;y_1,w)T(w)e_{v_0}(t-s;w,y_2)dw \ \ ,$$
which is exactly the kernel $\Ee_{v_0}$ .
\end{proof}

\begin{corollary}\label{c:1}
The spectral flow formula (\ref{e:3}) follows from the {\it EAF}.
\end{corollary}

\begin{proof} The
{\it EAF} shows that 
$k_{\e}(t;(v_1,y_1),(v_2,y_2))$ is equal to
$$2\e\alpha'(\e v_1)e_{\partial_v}(t;v_1,v_2)
(e_{v_1}*Te_{v_1})(t;y_1,y_2) + 
O(\frac {\e^{\frac 32}}{\sqrt{t}}) \ \ .$$
It follows that, for $\e$ small enough, 
we have the equality
$$index \ \Dd = 2\int_{N_{\e}}tr \ \e\alpha'(\e v_0)
e_{\partial_v}(t;v_0,v_0)
(e_{v_0}*Te_{v_0})(t;y,y)dydv_0$$
$$=\int_0^{\frac 1{\e}}
\frac {\e\alpha'(\e v_0)}{\sqrt{\pi t}}dv_0
\int_0^tTr_Me^{-sB_{\e v_0}^2}
Te^{-(t-s)B_{\e v_0}^2}ds$$ 
$$=\int_0^{\frac 1{\e}}
\frac {\e\alpha'(\e v_0)}{\sqrt{\pi t}}dv_0
\int_0^tTr_MTe^{-tB_{\e v_0}^2}ds$$
$$=\sqrt{\frac t{\pi}}\int_0^{\frac 1{\e}}\alpha'(\e v_0)
Tr_MTe^{-tB_{\e v_0}^2}(\e dv_0) = 
\sqrt{\frac t{\pi}}\int_0^1
Tr_M{\dot B}_ue^{-tB_u^2}du \ \ .$$
\end{proof}

Let us observe that 
(for compatible Dirac operators)
$$\lim_{t \to 0}\sqrt{\frac t{\pi}}{\cdot}
\int_0^1Tr_M{\dot B}_ue^{-tB^2_u}du = 
\int_0^1du\lim_{t \to 0}\sqrt{\frac t{\pi}}
{\cdot}Tr_M{\dot B}_ue^{-tB^2_u} = 
\int_0^1{\dot \eta}_udu \ \ ,$$
where $\eta_u = \frac 12 (dim \ ker(B_u) + \eta_{B_u}(0))$ 
is the $\eta$-invariant of 
the operator $B_u$ \cite{APS3}.

\subsection{Technicalities.}

We write the operator $\Dd_{\e}$ as
\begin{equation}\label{e:8}
\Dd_{\e} = \partial_v + (B_0 + \alpha(\e v_0)T) + 
(\alpha(\e v) - \alpha(\e v_0))T = \partial_v + B_{\e v_0} 
+ \beta (\e v)T
\ \, .
\end{equation}
In the following we consider 
the operator $\Dd_{\e}$ as 
an operator living on $\R \times M$ . 
This does not change anything in the proof, 
but simplifies the computations. 
Of course there is a problem with the definition 
of the index of $\Dd_{\e}$ in this set-up. 
Even though $B_{\e v}$ is a constant 
operator for $v < 0$ and $\frac 1{\e} \le v$ 
the index may not be well-defined unless 
the operator $B_{\e v}$ is invertible 
for those values of the normal coordinate. 
Hence, one can think that we perturbed 
the tangential operator by a small number and 
the invertibility condition is satisfied. 
In any case it follows that the integral from 
the kernel $k_{\e}(t;(v,y),(v,y))$ 
over $[0, \frac 1{\e}] \times M$ gives 
an integer equal to the index of $\Dd_{\e}$ 
on $N_{\e}$ and the integral over the leftover 
of $\R \times M$ gives a finite error term, which 
goes to $0$ as $\e \to 0$. Hence it is not difficult 
to show that
$$index \ \Dd_{\e} = 
\int_{\R \times M} 
tr \ k_{\e}(t;(v,y),(v,y))dydv \ ,$$
where on the left side we have 
the operator on $N_{\e}$.

We work on $\R \times M$ and we have to show 
that {\it EAF} holds at any given point 
$(v_0,y)$. After reparametrization, 
we can assume that $v_0 = 0$ and our operator 
has the form
$$\Dd_{\e} = \partial_v + B_0 + 
\beta(\e v)T \ \ ,$$
where the cut-off function $\beta(v)$ 
satisfies $\beta(0) = 0$. 
The corresponding Laplacians are
$$\Delta_{1,\e} = \Dd_{\e}^*\Dd_{\e} =
-\partial_v^2 + B_0^2 + 
\beta (\e v)(B_0T + TB_0) + \beta^2(\e v)T^2 - 
\e\beta '(\e v)T$$
$$= -\partial_v^2 + B_0^2 + \beta(\e v)T_1 - 
\e\beta '(\e v)T \ \ ,$$
and
$$\Delta_{2,\e} = \Dd_{\e}\Dd_{\e}^* = 
-\partial_v^2 + B_0^2 + \beta(\e v)T_1 + 
\e\beta '(\e v)T \ \ ,$$
where $T_1$ denotes the $1$st order 
tangential operator 
$B_0T + TB_0 + \beta(\e v)T^2$. 

 To evaluate $k_{\e}(t;(0,y),(0,y))$
we apply Duhamel's Principle \cite{McKS,BW}:
$$e^{-t\Delta_{1,\e}} - e^{-t\Delta_{2,\e}} = 
\int_0^t \frac d{ds}(e^{-s\Delta_{1,\e}}
e^{-(t-s)\Delta_{2,\e}})ds$$
$$= \int_0^t e^{-s\Delta_{1,\e}}(\Delta_{2,\e} - 
\Delta_{1,\e})e^{-(t-s)\Delta_{2,\e}}ds \ \ .$$
The difference $\Delta_{2,\e} - \Delta_{1,\e}$ 
is equal to the bundle endomorphism $2\e\beta'(\e v)T$ 
and this is the term which brings the $1$st, and the most important, 
factor of $\e$ into the formula. Once again, we simplify 
the presentation and introduce Laplacian 
$\Delta_0 = -\partial_v^2 + B_0^2$ and study each summand 
in the equality
$$e^{-t\Delta_{1,\e}} - e^{-t\Delta_{2,\e}} = 
(e^{-t\Delta_{1,\e}} - e^{-t\Delta_0}) + 
(e^{-t\Delta_0} - e^{-t\Delta_{2,\e}}) \ \ ,$$

The application of Duhamel's Principle to the first 
summand leads to the series
$$e^{-t\Delta_{1,\e}} - e^{-t\Delta_0} = 
\int_0^t e^{-s\Delta_{1,\e}}(\Delta_0 - 
\Delta_{1,\e})e^{-(t-s)\Delta_0}ds$$
$$\int_0^t e^{-s_1\Delta_0}(\Delta_0 - 
\Delta_{1,\e})e^{-(t-s_1)\Delta_0}ds_1$$ 
$$+ \int_0^tds_1\int_0^{s_1}ds_2 \ e^{-s_2\Delta_{1,\e}}
(\Delta_0 - \Delta_{1,\e})e^{-(s_1-s_2)\Delta_0}  
(\Delta_0 - \Delta_{1,\e})e^{-(t-s_1)\Delta_0}$$
$$= \sum_{k=1}^{\infty}\int_0^tds_1...\int_0^{s_{k-1}}ds_k \ 
e^{-s_k\Delta_0}(\Delta_0 - \Delta_{1,\e})...(\Delta_0 - 
\Delta_{1,\e})e^{-(t-s_1)\Delta_0} \ .$$

The result we state now is a standard application of 
Duhamel's Principle

\begin{prop}\label{p:1}
Let us consider the operator 
$\Delta_R = \Delta_0 + \beta(\e v)R$ , 
where $R : C^{\infty}(M;S) \to C^{\infty}(M;S)$ 
is a tangential differential operator of order $1$ . 
Then, there exists positive constants $c_1$ 
and $c_2$ (independent of $\e$) such that for 
any sufficiently small $\e$ the following 
estimate holds
\begin{equation}\label{e:9}
\|e_R(t;(v_1,y_1),(v_2,y_2))\| \le 
c_1t^{-\frac {m+1}2}e^{-c_2 
\frac {d^2(((v_1,y_1),(v_2,y_2))}2}
\ \, ,
\end{equation}
where $e_R(t;(v_1,y_1),(v_2,y_2))$ denotes 
the kernel of $e^{-t\Delta_R}$ .
\end{prop}

\begin{proof}
The estimate (\ref{e:9}) holds for $\Delta_0$ and we have 
$$e^{-t\Delta_R} - e^{-t\Delta_0} = 
\sum_{k=1}^{\infty}\int_0^tds_1...\int_0^{s_{k-1}}ds_k \ 
e^{-s_k\Delta_0}(\beta(\e v)R)...(\beta(\e v)R)e^{-(t-s_1)\Delta_0}$$
$$= \sum_{k=1}^{\infty}\int_0^tds_1...\int_0^{s_{k-1}}ds_k 
\biggm(e^{-s_k(-\partial_v^2)}(\beta(\e v)).....
(\beta(\e v))e^{-(t-s_1)(-\partial_v^2)}\biggm) \times $$
$$\biggm(e^{-s_kB_0^2}R....Re^{-(t-s_1)B_0^2}\biggm) \ \ .$$
The kernel of the operator in the first bracket is estimated as follows
$$\|\biggm(e^{-s_k(-\partial_v^2)}(\beta(\e v)).....
(\beta(\e v))e^{-(t-s_1)(-\partial_v^2)}\biggm)(t;v_1,v_2)\| \le $$
$$\|\biggm(e^{-s_k(-\partial_v^2)}...
e^{-(t-s_1)(-\partial_v^2)}\biggm)(t;v_1,v_2)\| \ ,$$
as $0 \le \beta(\e v)$. This leads to the estimate
$$\|e_R(t;(v_1,y_1),(v_2,y_2)) - e_0(t;(v_1,y_1),(v_2,y_2))\| \le$$ 
$$\frac {e^{- \frac {(v_1 - v_2)^2}{4t}}}{\sqrt{4 \pi t}}
\sum_{k=1}^{\infty} \| \biggm(e^{-s_kB_0^2}R....
Re^{-(t-s_1)B_0^2}\biggm)(t;y_1,y_2)\| \ .$$
The series here is estimated in the standard way. The kernel of 
the operator $Re^{-tB_0^2}$ is bounded by 
$c_1t^{-\frac {m+1}2}e^{- c_2\frac {(v_1 - v_2)^2}{t}}$. 
Therefore we follow \cite{McKS} and obtain
$$\sum_{k=1}^{\infty} \| \biggm(e^{-s_kB_0^2}R....
Re^{-(t-s_1)B_0^2}\biggm)(t;y_1,y_2)\| \le 
c_3t^{-\frac {m-1}2}e^{- c_4\frac {d^2(y_1,y_2)}{t}} \ .$$
The positive constants above  do not depend on $\e$.
\end{proof}

It is no problem to see that the estimate (\ref{e:9}) from 
Proposition \ref{p:1} is also satisfied by kernels of the heat 
operators defined by $\Del_{1,\e}$ and $\Del_{2,\e}$, which 
leads to the following useful property.

\begin{corollary}\label{c:2}
The contribution to 
the kernel $k_{\e}(t;(0,y_1),(0,y_2))$ 
provided by the points distant
more than $\frac 1{\sqrt{\e}}$ from $\{0\} \times M$ may be disregarded.
\end{corollary}

\begin{proof}
We have
$$\|\left(e^{-t\Del_{1,\e}} - 
e^{-t\Del_{2,\e}}\right)(t;(0,y_1),(0,y_2))\| =$$
$$\|\int_0^te^{-s\Del_{1,\e}}(2\e\beta(\e v)T)
e^{-(t-s)\Del_{2,\e}}(t;(0,y_1),(0,y_2))ds\| \le$$
$$\int_0^tds \int_{\R \times M}du dz 
\|e_{1,\e}(s;(0,y_1),(u,z))
2\e\beta(\e u)T(u)e_{2,\e}(t-s;(u,z),(0,y))\| \ \ .$$
We want to show that the integral over $|u| > \frac 1{\sqrt{\e}}$ 
is exponentially small with respect to $\e$.
$$\int_0^tds \int_{|u| > \frac1{\sqrt{\e}}}du \int_Mdz 
\|e_{1,\e}(s;(0,y_1),(u,z))\|
2\e\beta(\e u)\times$$
$$\times \|T(u)\|\|e_{2,\e}(t-s;(u,z),(0,y))\|$$
$$\le c_3\e\int_0^tds\int_{|u| > \frac 1{\sqrt{\e}}}du 
\int_M dz \ c_1s^{- \frac {m+1}2}
e^{-c_2\frac {d^2((0,y_1),(u,z))}s}
(t-s)^{- \frac {m+1}2}\times $$
$$\times e^{-c_2
\frac {d^2((u,z),(0,y_2))}{t-s}}$$
$$\le c_4\e{\cdot}vol(M)\int_0^tds
\int_{|u| > \frac 1{\sqrt{\e}}}(s(t-s))^{-\frac {m+1}2}
e^{-c_2\frac {u^2}s}e^{-c_2\frac {u^2}{t-s}}du $$
$$\le c_4\e{\cdot}vol(M)\int_0^tds\int_{|u| > \frac 1{\sqrt{\e}}}
e^{-c_5\frac {u^2}s}e^{-c_5\frac {u^2}{t-s}}du $$
$$\le c_4\e{\cdot}vol(M)
\int_0^tds\int_{|u| > \frac 1{\sqrt{\e}}}
e^{-c_5\frac {u^2}t}du \le c_6\e 
{\cdot}e^{- \frac {c_7}{\e t}} \ \ .$$
\end{proof}

Now the idea of the proof of the {\it EAF} 
can be easily understood. The kernel $k_{\e}(t;(0,y_1),(0,y_2)$ 
is the kernel of the operator
\begin{equation}\label{e:10}
\int_0^t e^{-s\Delta_{1,\e}}(\Delta_{2,\e} - 
\Delta_{1,\e})e^{-(t-s)\Delta_{2,\e}}ds
\ \, .
\end{equation}
Kernels of both heat operators $e^{-t\Delta_{i,\e}}$ , 
expands into series where the leading term is the kernel 
of the operator $e^{-t\Delta_0}$ . It follows 
that the leading term in the expansion 
of the kernel of the operator (\ref{e:10}) is
\begin{equation}\label{e:11}
\int_0^tds\int_{\R \times M}e_0(s;(0,y_1),(z,w))
(2\e \beta'(\e z)T(w))e_0(t-s;(z,w),(0,y_2))dzdw
\ \, .
\end{equation}
This gives us the contribution which appears in the {\it EAF}.

We have to show that further perturbation brings the consecutive 
powers of $\sqrt{\e}$ into the picture in order to finish the proof. 
So let us replace $e_0(s;(0,y_1),(z,w))$ 
by $e_{1,\e}(s;(0,y_1),(z,w))$ in the formula (\ref{e:11}), 
hence
$$\int_0^t e^{-s\Delta_0}(2\e \beta'(\e z)T(w))e^{-(t-s)\Delta_0}ds$$ 
is replaced by 
$$\int_0^t e^{-s\Delta_{1,\e}}(2\e \beta'(\e z)T(w))e^{-(t-s)\Delta_0}ds$$
and the kernel of $e^{-s\Delta_{1,\e}} - e^{-s\Delta_0}$
brings the extra factor $\sqrt{\e}$ as the elementary estimates 
presented below show. 
The operator $\Delta_{1,\e}$ is obtained from 
$\Delta_0$ by adding the correction term of the form 
$\g(\e v)S_1$ , where $\g(\e v)$ is the cut-off function with the properties 
specified earlier.

\bigskip

\noindent{\bf The main estimate}

We have
$$\|e_{1,\e}(t;(0,y),(0,w)) - e_0(t;(0,y),(0,q))\|$$
$$=\|\int_0^tds \int_{|v| \le \sqrt{\e}}dv \int_Mdz \ 
e_{1,\e}(s;(0,y),(v,z))\gamma(\e v)S_1(z)
e_0(t;(v,z),(0,q))\|$$
$$\le \int_0^tds \int_{-\frac 1{\sqrt{\e}}}^{-\frac 1{\sqrt{\e}}}dv 
\int_Mdz \ 
\|e_{1,\e}(s;(0,y),(v,z))\|{\cdot}|\gamma(\e v)|
{\cdot}\|S_1(z)e_0(t;(v,z),(0,q))\|$$
$$\le c\sqrt{\e}\int_0^tds 
\int_{-\frac 1{\sqrt{\e}}}^{-\frac 1{\sqrt{\e}}}dv \int_Mdz \ 
c_1s^{-\frac {m+1}2}
e^{-c_2 \frac {v^2 + d^2(y,z)}s}c_1
(t-s)^{-\frac {m+2}2}
e^{-c_2 \frac {v^2 + d^2(z,q)}{t-s}}$$
$$\le c\sqrt{\e}c_1^2\int_0^t\frac 1{\sqrt{t-s}}ds 
\int_{-\frac 1{\sqrt{\e}}}^{-\frac 1{\sqrt{\e}}}dv 
{\frac {e^{-c_2 \frac {tv^2}{s(t-s)}}}
{\sqrt{s(t-s)}}}
\int_Mdz \frac{e^{-c_2\frac {d^2(y,z)}s}
e^{-c_2\frac {d^2(z,q)}{t-s}}}{(s(t-s))^{-\frac m2}}$$
$$\le c_3 \frac {\sqrt{\e}}{\sqrt{t}}{\cdot}
\int_0^t\frac 1{\sqrt{t-s}}ds 
\int_Mdz \ (s(t-s))^{-\frac m2} 
e^{-c_2\frac {d^2(y,z)}s}
e^{-c_2\frac {d^2(z,q)}{t-s}} \ \ .$$

The following elementary inequality 
is used to estimate the factor 
$e^{-c_2\frac {d^2(y,z)}s}
e^{-c_2\frac {d^2(z,q)}{t-s}}$,
$$\frac {d^2(y,q)}t \le 
\frac {d^2(y,z)}s + \frac {d^2(z,q)}{t-s} \ \ .$$
We have
$$e^{-c_2\frac {d^2(y,z)}s}
e^{-c_2\frac {d^2(z,q)}{t-s}} = 
e^{-2c_2\frac {d^2(y,z)}{2s}}
e^{-2c_2\frac {d^2(z,q)}{2(t-s)}} \le$$
$$e^{-c_2\frac {d^2(y,q)}{2t}}
e^{-c_2\frac {d^2(y,z)}{2s}}
e^{-c_2\frac {d^2(z,q)}{2(t-s)}} \ \ .$$
All this amounts to
$$\|e_R(t;(0,y),(0,w)) - 
e_0(t;(0,y),(0,q))\| $$
$$\le c_3\sqrt{\frac {\e}{t}}{\cdot}
e^{-c_2\frac {d^2(y,q)}{2t}}
\int_0^t \frac {ds}{\sqrt{t-s}}
\int_Mdz \ (s(t-s))^{-\frac m2} 
e^{-c_2\frac {d^2(y,z)}{2s}}
e^{-c_2\frac {d^2(z,q)}{2(t-s)}}$$ 
$$\le c_4\sqrt{e}{\cdot}
e^{-c_2\frac {d^2(y,q)}{2t}}
\int_Mdz \ (s(t-s))^{-\frac m2} 
e^{-c_2\frac {d^2(y,z)}{2s}}
e^{-c_2\frac {d^2(z,q)}{2(t-s)}}$$
$$ 
\le c_5\sqrt{e}
{\cdot}t^{-\frac m2}
{\cdot}e^{-c_2\frac {d^2(y,q)}{2t}} \ \ .$$

Now, we can show the fact used in the proof of the {\it EAF} 
in the previous subsection.
We have
$$e^{-t\Del_{1,{\e}}} - e^{-t\Del_{2,{\e}}} = 
\int_0^t \ e^{-s\Del_{1,{\e}}}
(\Del_{2,{\e}} - \Del_{1,{\e}})
e^{-(t-s)\Del_{2,{\e}}}ds$$
$$=\int_0^t [e^{-s\Del_0}(\Del_{2,{\e}} - \Del_{1,{\e}})
e^{-(t-s)\Del_{2,{\e}}}+ 
(e^{-s\Del_{1,\e}} - e^{-s\Del_0})
(\Del_{2,{\e}} - \Del_{1,{\e}})
e^{-(t-s)\Del_{2,{\e}}}]ds$$
$$=\int_0^te^{-s\Del_0}(\Del_{2,{\e}} - \Del_{1,{\e}})
e^{-(t-s)\Del_0}ds $$ $$+ 
\int_0^t (e^{-s\Del_{1,\e}} - e^{-s\Del_0})
(\Del_{2,{\e}} - \Del_{1,{\e}})
e^{-(t-s)\Del_{2,{\e}}}ds $$
$$+ \int_0^te^{-s\Del_0}
(\Del_{2,{\e}} - \Del_{1,{\e}})
(e^{-(t-s)\Del_{2,{\e}}} - 
e^{-(t-s)\Del_0})ds \ \ .$$

The operator 
$\Del_{2,{\e}} - \Del_{1,{\e}} = 2\e \beta'(\e v)T$ 
brings a factor $\e$ and the first term 
on the right side is of the form 
$\int_0^te^{-s\Del_0}(\Del_{2,{\e}} - 
\Del_{1,{\e}})e^{-(t-s)\Del_0}ds$ 
and of the $\e$ size. The next term 
on the right side 
$$\int_0^t (e^{-s\Del_{1,\e}} - e^{-s\Del_0})
(\Del_{2,{\e}} - \Del_{1,{\e}})e^{-(t-s)\Del_{2,{\e}}}ds$$
contains additionally the difference 
$e^{-s\Del_{1,\e}} - e^{-s\Del_0}$, 
hence it is of the size $\e^{\frac 32}$. 
This is also the case of the last term 
$\int_0^te^{-s\Del_0}(\Del_{2,{\e}} - \Del_{1,{\e}})
(e^{-(t-s)\Del_{2,{\e}}} - e^{-(t-s)\Del_0})ds$. 
It follows that, as we take limit as $ \e \to 0$, 
only the integral (over $\R \times M$) from 
the first term is going to survive. 
The kernel of this operator 
at the point $(t,(0,y_1),(0,y_2))$ 
has the form
$$(e_0\#(2\e \beta'(\e v)Te_0)
(t,(0,y_1),(0,y_2)) \ \ ,$$
where $e_0(t;(v_1,y_1),(v_2,y_2))$ 
denotes kernel of the operator $e^{-tB^2_0}$. 
This is exactly what we need in order to complete the proof.

\subsection{The EAF for operators on covering spaces}

We let $\tilde M$ be the universal covering space 
for the closed manifold $M$~, with the corresponding
fundamental group $G (= \pi_1(M))$. We assume that we are given  
a $G$-invariant, compatible Dirac operator on $\tilde M$
$${\tilde B} : C^{\infty}({\tilde M};{\tilde S}\otimes {\tilde E}) 
\to C^{\infty}({\tilde M};{\tilde S}\otimes {\tilde E}) \ \ ,$$
where $G$ acts on $\tilde E$ via a representation $\rho$~.
The appropriate von Neumann algebra $\maN$ is the commutant of
the $G$ action and there is a corresponding
trace $\tau= \tau_G$
as described by Atiyah \cite{A} and introduced at the beginning of this Section.


We introduce now a $G$-invariant unitary bundle automorphism 
$\tilde h$ of the auxiliary bundle $\tilde E$~, and we consider 
family $\{{\tilde B}_u\}$ as defined in (\ref{e:g1}). This family 
has a well-defined spectral flow and we want to prove that
\begin{equation}\label{e:51}
sf\{ \tilde B_s\} = index\tilde A
\end{equation}
(Here the tilde just denotes the covering space 
analogues of the operators we  introduced before).
This is the result which corresponds 
to Theorem \ref{t:1}
\begin{theorem}\label{t:51}
\begin{equation}\label{e:52}
index\tilde A = sf\{ \tilde B_s\} = \int_0^1{\dot \eta}^G_sds.
\end{equation}
\end{theorem}
We remind the reader that the $\eta$-invariant in this context 
was studied by Cheeger and Gromov \cite{ChG},
Hurder \cite{H} and later on by Mathai \cite{M} and others.
To prove the theorem \ref{t:51} have only to show 
\begin{equation}\label{e:eafc1}
index_G\tilde A =\int_0^1\sqrt{\frac t{\pi}}
{\cdot}\tau_G{\dot B}_s
e^{-tB_s^2}ds,
\end{equation}
This however follows easily from the extension of the {\it EAF} 
to the present context. 
We follow the previous argument.
The $index_G{\tilde A}$  is equal to
$$\tau_G(e^{-t{\tilde A}^*{\tilde A}} - 
e^{-t{\tilde A}{\tilde A}^*}) $$
(The McKean-Singer formula 
for the index 
holds when the von Neumann algebra is not necessarily a factor \cite{CPRS3}). 
We blow up the metric and all the arguments 
from the compact case come through with the slight 
modifications. 

The only problem we face is that although our trace is defined 
by the integral over the (compact) fundamental domain we have to integrate 
over the whole non-compact manifold $\tilde M$~, when applying 
{\it Duhamel's Principle}. 
First let us notice that the standard point-wise estimate on 
the heat-kernel (\ref{e:9}) holds on $\tilde M$ \cite{D79}. 
The difficulty follows from the well-known fact that the volume of 
the ball with a fixed center on $\tilde M$ may grow exponentially 
with the radius of the ball. Therefore we have to be careful with 
the arguments which lead to the proof of the results which correspond
to Proposition \ref{p:1} and Corollary \ref{c:2}. 
Actually, everything 
works only because the volume growth is at most exponential
in the diameter.
We omit the estimates 
which lead to the proof of Proposition \ref{p:1} for the case of the
covering space $\tilde M$. However to present the flavour of the computations 
involved we present a modification of the argument used to get 
Corollary \ref{c:2}. We work on the space ${\tilde M}_{\e}$ (or on 
$\R \times {\tilde M}$) now, so 
the proper formulation of the result is as follows.

\begin{corollary}\label{c:3}
Let $y_1 , y_2 \in {\bar M}$ then we may disregard the contribution to 
the kernel $k_{\e}(t;(0,y_1),(0,y_2))$ provided by the points 
$(u,z) \in \R \times {\tilde M}$ such that $|u| > \frac 1{\sqrt{\e}}$ and 
$d(z,{\bar M}) > 1$.
\end{corollary}

\begin{proof}
We start as in Section 8.2
$$\|e^{-t\Del_{1,\e}} - 
e^{-t\Del_{2,\e}}\|(t;(0,y_1),(0,y_2)) =$$
$$\|\int_0^te^{-s\Del_{1,\e}}(2\e\beta(\e v)T)
e^{-(t-s)\Del_{2,\e}}\|(t;(0,y_1),(0,y_2))ds \le$$
$$\int_0^tds \int_{\R \times {\tilde M}}du dz 
\|e_{1,\e}(s;(0,y_1),(u,z))
2\e\beta(\e u)T(u)e_{2,\e}(t-s;(u,z),(0,y))\|.$$

We have to remember that nothing is changed (i.e. we obtain a negligible 
contribution), while we work on compact pieces of the space $\R \times 
{\tilde M}$~. 
Here we want to show that the integral over 
$(u,z) \in \R \times {\tilde M}$ such that $|u| > \frac 1{\sqrt{\e}}$ and 
$d(z,{\bar M}) > 1$ 
is exponentially small with respect to $\e$ .
$$\int_0^tds \int_{|u| > \frac1{\sqrt{\e}}}du 
\int_{\{z ; d(z,{\bar M}) >1\}}dz$$ 
$$\|e_{1,\e}(s;(0,y_1),(u,z))\|
2\e\beta(\e u)\|T(u)\|\|e_{2,\e}(t-s;(u,z),(0,y))\| \le$$
$$c_3\e{\cdot}\int_0^tds\int_{|u| > \frac 1{\sqrt{\e}}}du 
\int_{\{z ; d(z,{\bar M}) >1\}} dz 
e^{-c_2\frac {d^2((0,y_1),(u,z))}s}e^{-c_2\frac {d^2((u,z),(0,y_2))}{t-s}} \le$$
$$c_3\e{\cdot}\int_0^tds
\int_{|u| > \frac 1{\sqrt{\e}}}
e^{-c_2\frac {u^2}s}e^{-c_2\frac {u^2}{t-s}}du \int_{\{z ; d(z,{\bar M}) >1\}} 
dz \ 
e^{-c_2\frac {d^2(z;{\bar M}}s}e^{-c_2\frac {d^2(z;{\bar M}}{t-s}}.$$

The first integral on the right side is estimated as in Section 2
$$\int_{|u| > \frac 1{\sqrt{\e}}}e^{-c_2\frac {u^2}s}e^{-c_2\frac {u^2}{t-s}}du 
\le$$
$$\int_{|u| > \frac 1{\sqrt{\e}}}e^{-c_2\frac {u^2}t}du \le 
\sqrt{\frac t{c_2}}{\cdot}\int_{|v| > \sqrt{\frac {c_2}t}}e^{-v^2}dv \le 
\sqrt{\frac t{c_2}}{\cdot}e^{-\frac {c_2}{\e t}}  .$$

The second integral involves volume of the manifold $\tilde M$. Modulo 
negligible 
error (up to a contribution from a compact set)
we can look at it as the integral over the outside of the ball centered at 
the fixed point ${\bar y} \in {\bar M}$ 
with radius $R = 1 + diam \ {\bar M}$ . We do have 
$$\int_{\{z ; d(z,{\bar y}) >R\}}  e^{-c_2\frac {d^2(z;{\bar y}}s}
e^{-c_2\frac {d^2(z;{\bar y}}{t-s}}dz \le \int_{\{z ; d(z,{\bar y}) >R\}}  
e^{-c_2\frac {d^2(z;{\bar y}}t}dz \le$$
$$c_3\int_R^{\infty}e^{-c_2\frac {r^2}t}e^{c_4r}dr \le 
c_5\int_R^{\infty}e^{-c_6\frac {r^2}t}dr \le$$
$$c_5\sqrt{\frac t{c_6}}{\cdot}\int_{v>\sqrt{\frac {c_6}t}R}e^{-v^2}dv \le 
c_5\sqrt{\frac t{c_6}}{\cdot}e^{-c_6\frac {R^2}t}$$
\end{proof}

\section{Spectral flow for almost periodic gauge transformations}

\subsection{Shubin's framework}

We follow Shubin \cite{Shu,Shu2} which in turn extends 
the original paper of Coburn et al \cite{CDSS,CMS}.
In this paragraph, we review the definition of the von Neumann algebra  which 
is appropriate for the study of almost periodic operators. Recall that a 
trigonometric function is a finite linear combination of exponential 
functions $e_\xi:x\mapsto e^{i<x,\xi>}$. The space $\Trig (\R^n)$ of 
trigonometric functions is clearly a $*-$subalgebra of the $C^*-$algebra 
$C_b(\R^n)$ of continuous bounded functions. The uniform closure of  
$\Trig (\R^n)$ is thus a $C^*-$algebra called the algebra of almost periodic 
functions and denoted $\CAP(\R^n)$. Since this $C^*-$algebra is unital and 
commutative, it is the $C^*-$algebra of continuous functions on a compact 
space $\R^n_B$ which is a compactification of $\R^n$ with respect to the 
appropriate topology. The compact space $\R^n_B$ is called the Bohr 
compactification of $\R^n$ or simply the Bohr space. Addition extends to 
$\R^n_B$ which is a compact abelian group containing $\R^n$ as a dense 
subgroup. There is a unique normalized Haar measure $\alpha_B$ on $\R^n_B$ 
such that the family $(e_\xi)_{\xi\in \R^n}$ is orthonormal. Namely, the 
measure $\alpha_B$ is given for any almost periodic function $f$ on $\R^n$ by:
$$
\alpha_B (f) := \lim_{T\to +\infty} \frac{1}{(2T)^n} \int_{(-T,T)^n} f(x) dx.
$$
By using the measure $\alpha_B$ one defines the Hilbert space completion 
$L^2(\R^n_B)$ of $\Trig (\R^n)$. This Hilbert space is called the Besicovich 
space and it has an orthonormal basis given by $(e_\xi)_{\xi\in \R^n}$. In 
other words, the Pontryagin dual of $\R^n_B$ is the discrete abelian group 
$\R^n_d$ and the Fourier transform $\maF_B:\ell^2(\R^n_d) \longrightarrow 
L^2(\R^n_B)$ is given by:
$$
\maF_B (\delta_\xi) = e_\xi, \quad \text{ with }\delta_\xi (\eta) = 
\delta_{\xi, \eta},
$$
where $\delta_{\xi, \eta}$ is the Kronecker symbol. We shall denote by $\maF$ 
the usual Fourier transform on the abelian group $\R^n$ with its usual 
Lebesgue measure.

For any $f\in C_b(\R^n)$ we shall denote, for any vector $\lambda\in \R^n$, 
by $T_\lambda f$ the translated function defined by 
$(T_\lambda f) (x) = f(x-\lambda)$. Let $f\in C_b(\R^n)$ and let $\ep >0$ be 
given. A vector $\lambda\in \R^n$ is called  an $\ep-$ period for $f$ if the 
uniform norm of $T_\lambda f - f$ is bounded by $\ep$, i.e.
$$
\|T_\lambda f - f\|_\infty := \sup_{t\in \R^n} |f(t-\lambda) - f(t)| \leq \ep.
$$
A subset $E$ of $\R^n$ is relatively dense if there exists $T > 0$ such that
$$\forall x\in \R^n, \exists u\in E: u-x \in \left[-\frac{T}{2}, 
+\frac{T}{2}\right]^n.$$
It is worth pointing out that,  for any function $f: \R^n \to \C$, 
the following properties are equivalent \cite{CDSS}:
\begin{itemize}
\item $f$ is an almost periodic function.
\item $f$ is a continuous bounded function whose $\ep$ periods are relatively 
dense for every $\ep >0$. 
\end{itemize}

It is clear from the second characterization of an almost periodic function 
that any periodic function is almost periodic. An interesting class of 
examples arises from the study of quasi-periodic functions. Assume for 
simplicity that $n=1$ and let $\alpha=(\alpha_1, \cdots, \alpha_p)\in \R^p$ be 
a fixed list of real numbers. Then for any summable sequence 
$c=(c_m)_{m\in \Z^p}$, we get an almost periodic function on $\R$ by setting:
$$
\varrho_{c,\alpha} (x) := \sum_{m\in \Z^p} c_m e^{2i\pi <m,\alpha> x}.
$$
Then more complicated examples of almost periodic functions arise as limits of 
periodic or quasi-periodic functions. For instance, the function 
$\sum_{n\geq 0} a_n \cos (\frac{x}{2^n})$ where $\sum_n |a_n| < +\infty$, is 
an almost periodic function. 

The action of $\R^n$ on $\R^n_B$ by translations yields a topological 
dynamical system whose naturally associated  von Neumann algebra is the 
crossed product von Neumann algebra $L^\infty (\R^n_B)\rtimes \R^n$.
It is more convenient for applications to
consider the commutant of this von Neumann algebra
 denoting it by $\cnn$. It is also a  crossed product. 
This time it is the von Neumann algebra $L^\infty (\R^n)\rtimes \R^n_d$. 
The von Neumann algebra $\cnn$ is a type II$_\infty$ factor with a 
faithful normal semi-finite trace $\tau$. It can be described as the set of 
Borel essentially bounded families $(A_\mu)_{\mu\in\R^n_B}$ of bounded 
operators in $L^2(\R^n)$ which are $\R^n-$equivariant, i.e.
such that
$$
A_\mu = \sigma_\mu (A_0)=T_{-\mu} A_0T_{\mu}, \quad \forall \mu\in \R^n.
$$
Here and in the sequel we denote by $\sigma_\mu$ conjugation 
of any operator with the translation $T_\mu$ so that
$
\sigma_\mu (B) = T_{-\mu} B T_\mu.
$
If we denote by $M_\varphi$ the operator of multiplication by a bounded  
function $\varphi$, then examples of such families are given for any $\lambda$ 
by the families
$$
(\sigma_\mu(M_{e_\lambda}))_{\mu\in \R^n_B}. 
$$
We choose the Fourier transform 
$${\mathcal F}f(\zeta)=\frac{1}{(2\pi)^{n/2}}\int_{\R^n}e^{ix\zeta}f(x)dx.$$
Now  the von Neumann algebra
$\cnn$ can be defined \cite{CDSS,CMS} as the double commutant
of the set of operators $\{ M_{e_\lambda} \otimes
M_{e_\lambda} , T_\lambda \otimes 1 | \lambda\in\R^n\}$  
on the Hilbert space $\maH=L^2(\R^n)\otimes
L^2(\R^n_B)$ 

There is a natural way to imbed 
the $C^*-$algebra $\CAP (\R^n)$ in $\cnn$ by setting
$$
\pi(f) := (\sigma_\mu (M_f))_{\mu\in \R^n_B}
$$ 
This family then
belongs to $\cnn$ and $\pi$ is clearly faithful. 
Viewed as an operator on $\maH$, $\pi(f)$ is given by $\pi(f) (g) (x,\mu) = 
f(x+\mu) g(x,\mu)$.
If $B=(B_\mu)_\mu$ is a positive  element of $\cnn$, then  we define the 
expectation $E(B)$ as the Haar integral:
$$E(B) := \int_{\R^n_B} B_\mu d\alpha_B(\mu).$$
Since the family $B$ is translation equivariant and since $\alpha_B$ is 
translation invariant, the operator $E(B)$ clearly commutes with the 
translation in $L^2(\R^n)$ and is therefore given by a Fourier multiplier 
$\wM (\varphi_B)$ with $\varphi_B$ a positive element of $L^{\infty} (\R^n)$. 
Recall that the Fourier multiplier $\wM (\varphi_B)$ is conjugation of the 
multiplication operator $M_\varphi$ by the Fourier transform, i.e. $
\wM (\varphi_B) = \maF^{-1} M_\varphi \maF.$
When the function $\varphi$ is for instance in the Schwartz space, the 
operator $\wM (\varphi_B)$ is convolution by the Schwartz function 
$\frac{1}{(2\pi)^{n/2}}\maF^{-1} \varphi$. 
Hence the expectation $E$ takes values in the von Neumann algebra 
$\wM (L^{\infty}(\R^n))$, i.e.
$$E: \cnn \longrightarrow \wM (L^{\infty}(\R^n)).$$
Now, using the usual
Lebesgue integral on $\R^n$, we use the normalisation
of Coburn et al \cite{CDSS}
and introduce the following definition of the trace $\tau$:
$$
\tau (B) = \int_{\R^n} \varphi_B (\zeta) d\zeta.
$$

\begin{lemma} \cite{CDSS,CMS}\
The map $\tau$ is, up to constant, the unique positive normal faithful 
semi-finite trace on $\cnn$.
\end{lemma}

The space $L^1(\cnn, \tau)$ of trace-class $\tau-$measurable operators with 
respect to $\cnn$ is the space of $\tau-$measurable operators $T$
as explained by Fack et al \cite{FK} 
such that $\int_{(0, +\infty)} \mu^\tau_s(T) ds < +\infty$. Here 
$\mu^\tau_s(T)$ is the $s-$th characteristic value of $T$ \cite{FK}, for the 
precise definitions. More generally and for any $p\geq 1$, we shall denote by 
$L^p(\cnn, \tau)$ the space of $\tau-$measurable operators $T$ such that 
$(T^*T)^{p/2} \in L^1(\cnn,\tau)$. It is well known that the space 
$L^p(\cnn,\tau)\cap \cnn$ 
is a two-sided $*-$ideal in $\cnn$ that we shall call 
the $p-$th Schatten ideal of $\cnn$. 

We also consider the Dixmier space $L^{1,\infty}(\cnn, \tau)$ of those 
operators $T\in \cnn$ such that
$$
\int_0^s \mu^\tau_t(T) dt \sim O(\log (s)).
$$
Again, $L^{1,\infty}(\cnn, \tau)$ is a two-sided $*-$ideal in $\cnn$. There 
are well defined Dixmier traces $\tau_\omega$ on $L^{1,\infty}(\cnn, \tau)$ 
parametrized by limiting processes $\omega$ \cite{BeF,CPS2}.

Consider the trace on the von Neumann algebra $\cnn$ evaluated on an
operator of the form $M_aK$ where $a$ is almost periodic and
$K$ is a convolution operator on 
$L^2(\R^n)$ arising from multiplication by an $L^1$ function
$k$ on the Fourier transform.
We have,
$$\tau(M_aK)
=  \lim_{T\to +\infty} \frac{1}{(2T)^n} \int_{(-T, T)^n}
a(x) dx\int_{\R^n}k(\zeta)d\zeta$$

More generally, any pseudodifferential operator $A$
on $L^2(R^n,\C^N)$ 
with almost periodic coefficients of nonpositive order $m$ acting on 
$\C^N-$valued functions, can be viewed as a family over $\R_{B}^n$ of 
pseudodifferential operators on $\R^n$. To do this first take
 the symbol  $a$ of $A$, then the operator $\sigma_\mu(A)$ is the 
pseudodifferential operator with almost periodic coefficients whose symbol is 
$$
(x,\xi) \longmapsto a(x+\mu, \xi).
$$
When $m\leq 0$, we get in this way an element of the von Neumann algebra 
$\cnn$. We denote by $\Psi_{AP}^0$ the algebra of pseudodifferential operators 
with almost periodic coefficients and with non positive order. 
When the order $m$ of $A$ is $>0$ then the operator $A^\sharp$
given by the family $(\sigma_\mu(A))_{\mu\in \R^n_B}$ is affiliated with the 
von Neumann algebra $\cnn$. If the order $m$ of $A$ is $<-n$, then the bounded 
operator $A^\sharp$ is trace class with respect to the trace $\tau$ on the von 
Neumann algebra $\cnn\otimes M_N(\C)$ \cite{Shu3}[Proposition 3.3] and we  
have:
$$
\tau (A^\sharp) = 
\lim_{T\to +\infty}\frac{1}{(2T)^n}\int_{(-T,+T)^n\times \R^n} \tr (a(x,\zeta)) 
dx d\zeta$$
Indeed, the expectation $E(A^\sharp)$ is a pseudodifferential operator on 
$\R^n$ with symbol denoted by $E(a)$ and is independent of the $x-$variable, 
it is given by:
$$E(a) (\zeta) = \lim_{T\to +\infty} \frac{1}{(2T)^n} 
\int_{(-T,+T)^n} a(x,\zeta) dx.$$
Hence the operator $E(A^\sharp)$ is precisely the Fourier multiplier 
$\wM (E(a))$ and so:
$$
\tau (A^\sharp) = \int_{\R^n} \tr( E(a) (\zeta)) d\zeta.
$$
Let $\Psi^\infty_{AP}$ be the space of one step polyhomogeneous classical 
pseudodifferential operators on $\R^n$ with almost periodic coefficients.

\begin{theorem}\label{Dixmier}\
Let $A$ be a (scalar) pseudodifferential operator with almost periodic 
coefficients on $\R^n$. We assume that the order $m$ of $A$ is  $\leq -n$ and 
we denote by $a_{-n}$ the $-n$ homogeneous part of the symbol $a$. Then the 
operator $A^\sharp$ belongs to the Dixmier ideal $L^{1,\infty}(\cnn, \tau)$. 
Moreover, the Dixmier trace $\tau_\omega (A^\sharp)$  of $A^\sharp$ associated 
with a limiting process $\omega$ does not depend on $\omega$ and is given by 
the formula:
$$
\tau_\omega (A^\sharp) = \frac{1}{n} \int_{\R^n_B \times 
\IS^{n-1}}
a_{-n} (x, \zeta) d\alpha_B(x) d\zeta.
$$
\end{theorem}

\begin{proof}\
We denote as usual by $\Delta$ the Laplace operator on $\R^n$. The 
operator $A(1+\Delta)^{n/2}$ is then a pseudodifferential operator with almost 
periodic coefficients and nonpositive order. 
Hence, the operator $[A(1+\Delta)^{n/2}]^\sharp=
A^\sharp(1+\Delta^\sharp)^{n/2}$ belongs to the von Neumann algebra $\cnn$. 
Now the operator $(1+\Delta^\sharp)^{-n/2}$ is a Fourier multiplier defined 
by the function $\zeta\mapsto 
(1+\zeta^2)^{n/2}$. 
Hence if, for $\lambda >0$, $E_\lambda$ is the spectral projection 
of the operator $(1+\Delta)^{-n/2}$ corresponding to the interval 
$(0,\lambda)$ then the operator $1-E_\lambda$ is the Fourier multiplier 
defined by the function 
$\zeta \mapsto 1_{(\lambda,+\infty)}(\zeta^2+1)^{n/2})$. 
It follows that the trace $\tau$ of the operator 
$1-E_\lambda$ is given by
$$
\int_{\R^n} 1_{(\lambda,+\infty)}(\frac{1}{(\zeta^2+1)^{n/2}}) d\zeta.
$$
It is easy to compute this integral and to show that it is proportional to 
$\frac{1}{\lambda}$. So the infimum of those $\lambda$ for which 
$\tau (1-E_\lambda) \leq t$ is precisely proportional to $\frac{1}{t}$. Hence 
the operator $(1+\Delta^\sharp)^{-n/2}$, and hence $A$,  belongs to the 
Dixmier ideal $L^{1,\infty}(\cnn, \tau)$.

In order to compute the Dixmier trace of the operator $A$, we apply 
Shubin \cite{Shu}[Theorem 10.1] to deduce that the spectral 
$\tau-$density $N_A(\lambda)$ of $A$ has the asymptotic expansion
$$
N_A(\lambda) = \frac{\chi_0(A)}{\lambda} (1+ o(1)), \quad \lambda \to +\infty,
$$
where $\chi_0(A)$ is given by:
$$
\chi_0(A) = \frac{1}{n} \int_{\R^n_B\times \IS^{n-1}} a_{-n} (x, \zeta) 
d\alpha_B(x) d\zeta.
$$
Now, if $A$ is positive then by Benameur et al \cite{BeF}[Proposition 1]:
$$
\tau_\omega (A) = \lim_{\lambda\to +\infty} \lambda N_A(\lambda) = \chi_0(A).
$$
This proves the theorem for positive $A$. Since the principal symbol map is a 
homomorphism, we deduce the result for general $A$.
\end{proof}

The reader familiar with the Wodzicki residue will observe that the 
normalisation we have chosen for the trace in the von Neumann setting of 
this Section
eliminates a factor of $\frac{1}{(2\pi)^n}$ which occurs at the corresponding 
point in the type I theory.

\section{The odd semifinite local index theorem}

The original type I version of this result is due to Connes-Moscovici
 \cite{CM}. There are two new proofs, one due to Higson \cite{Hi} and one
due to Carey et al \cite{CPRS2}. 
The latter argument handles the case of semifinite
spectral triples. Quite remarkably this very general  
odd semifinite local index theorem is proved by starting from the integral 
formulae for spectral flow that we have described in earlier sections.
We do not have the space here to explain how it is done.

We restrict our discussion to a statement of the theorem.
First, we require multi-indices $(k_1,...,k_m)$, $k_i\in\{0,1,2,...\}$, whose
length $m$ will always be
clear from the context. We write $|k|=k_1+\cdots+k_m$, and define
$\alpha(k)$ by
$$ \alpha(k)=1/{k_1!k_2!\cdots
k_m!(k_1+1)(k_1+k_2+2)\cdots(|k|+m)}.$$
The numbers $\sigma_{n,j}$ are defined by the equality
$$\prod_{j=0}^{n-1}(z+j+1/2)=\sum_{j=0}^{n}z^j\sigma_{n,j}$$
with $\sigma_{0,0}=1$.
These are just the elementary symmetric functions of $1/2,3/2,...,n-1/2$.

If $(\maA,\HH,\D)$ is a smooth semifinite spectral triple
(ie $\maA$ is in the domain of $\delta^n$
for all $n$ where $\delta(a)= [(1+D^2)^{1/2}, a]$)
and $T\in\cnn$, we write
$T^{(n)}$ to denote  the iterated commutator
$$[\D^2,[\D^2,[\cdots,[\D^2,T]\cdots]]]$$
 where we have $n$ commutators with
$\D^2$. It follows \cite{CPRS2}
that operators of the form 
$$T_1^{(n_1)}\cdots T_k^{(n_k)}(1+\D^2)^{-(n_1+\cdots+n_k)/2}$$ 
are in $\cnn$ when
$T_i=[\D,a_i]$, or $=a_i$ for $a_i\in\maA$.

\begin{definition}\label{dimension}
If $(\maA,\HH,\D)$ is a smooth semifinite spectral triple,
we call
$$ p=\inf\{k\in{\R}:\tau((1+\D^2)^{-k/2})<\infty\}$$
the {\bf spectral dimension} of $(\maA,\HH,\D)$.
We say that $(\maA,\HH,\D)$ has {\bf isolated spectral dimension} if
for $b$ of the form
$$b=a_0[\D,a_1]^{(k_1)}\cdots[\D,a_m]^{(k_m)}(1+\D^2)^{-m/2-|k|}$$
the zeta functions
$$ \zeta_b(z-(1-p)/2)=\tau(b(1+\D^2)^{-z+(1-p)/2})$$
have analytic continuations to a deleted neighbourhood of $z=(1-p)/2$.
\end{definition}

Now we define, for $(\maA,\HH,\D)$ having isolated spectral dimension and
$$b=a_0[\D,a_1]^{(k_1)}\cdots[\D,a_m]^{(k_m)}(1+\D^2)^{-m/2-|k|}$$
$$ \tau_j(b)=res_{z=(1-p)/2}(z-(1-p)/2)^j\zeta_b(z-(1-p)/2).$$
The hypothesis of isolated spectral dimension is clearly necessary here in 
order to define the residues. The
semifinite local index theorem is as follows.
\begin{theorem}\label{SFLIT} 
Let $(\maA,\HH,\D)$ be an 
odd finitely summable smooth spectral triple with spectral
dimension $p\geq 1$. Let $N=[p/2]+1$ where $[\cdot]$ denotes the 
integer part (so $2N-1$ is the largest odd integer $\leq p+1$), 
and let $u\in\maA$ be unitary. Then
 if $(\maA,\HH,\D)$ also has isolated spectral dimension then
$$ sf(\D,u^*\D u)=\frac{1}{\sqrt{2\pi i}}\sum_m 
(-1)^{(m-1)/2}(\frac{(m-1)}{2})!\phi_m(u,u^*,\ldots,u, u^*)$$
where 
$\phi_m(u,u^*,\ldots,u, u^*)$
is 
$$
\sum_{|k|=0}^{2N-1-m}
\sum_{j=0}^{|k|+(m-1)/2}(-1)^{|k|}
\alpha(k)\sigma_{(|k|+(m-1)/2),j}$$
$$\times \tau_j
\left(u[\D,u^*]^{(k_
1)}\cdots[\D,u]^{(k_m)}(1+\D^2)^{-|k|-m/2}\right),$$
When $[p]=2n$ 
is even, the term with $m=2N-1$ is zero, and 
for $m=1,3,...,2N-3$, all the top terms with $|k|=2N-1-m$ are zero. 
\end{theorem}

We aim to compute the terms in this formula for semifinite spectral
flow in the case where $\cD$ is the Euclidean Dirac operator on
the spin bundle $\cS$ over $\R^n$ tensored with the trivial bundle
rank $N$ bundle
and $u$ is a smooth almost periodic
function from $\R^n$ to $U(N)$.

\section{Almost periodic spectral triple}

We now apply
 the local index theorem to compute
spectral flow. We thus assume that $n$ is odd.
The von Neumann algebra constructed previously is non-separable
and so to avoid a discussion of the non-separable situation
we need to slightly modify our approach in this Section.
In fact it is sufficient to study the
dense countable abelian subgroups of
$\R^n$. Let us fix one such, call it $D$
and explain how the theory works for this case.
Consider the subalgebra $\mathcal A$ of
$\CAP(\R^n)$ consisting of almost periodic functions
generated by $e_\xi$ with $\xi\in D$.
 We denote by $\cA^\infty$ the $*$-subalgebra
of $\CAP(\R^n)$
consisting of functions in  $\mathcal A$
which have bounded derivatives of all orders.
The von Neumann algebra we now consider
is the crossed product algebra of $D$ with
$L^\infty(\R^n)$)and is denoted by $\cM$.
We take the Hilbert space on which this algebra acts to be
$B_D^2(\R^n)\otimes L^2(\R^n)$ where $B_D^2(\R^n)$
is the Hilbert space completion of $\mathcal A$
 where the
norm and inner product are given by the restriction of the Haar
trace on $\CAP^{\infty}(\R^n)$ to $\cA$
(note that $B_D^2(\R^n)\cong \ell^2(D)$).
This type II$_{\infty}$ von Neumann algebra
is endowed with a faithful normal semi-finite trace that we denote by
$\tau$.
(We note that the explicit formula for $\tau$ is as given in
Section 9.)


The usual Dirac operator on $\R^n$ is denoted by $\eth$. So, if $\cS$ 
carries the spin representation of the Clifford algebra of
$\R^n$ then $\eth$ acts on smooth $\cS-$valued functions on $\R^n$. The 
operator $\eth$ is $\Z^n-$periodic and it is affiliated with the von Neumann 
algebra $\cM_\cS = \cM \otimes \End(\cS)$. This latter is also a 
type II$_\infty$ von Neumann algebra with the trace $\tau\otimes \tr$. More 
generally, for any $N\geq 1$, we shall denote by $\cM_{\cS, N}$ the von 
Neumann algebra $\cM\otimes \End(\cS\otimes \C^N)$ with the trace 
$\tau\otimes \tr$.

The algebra $\cA$ and its closure are faithfully represented 
as $*-$subalgebras of the von Neumann algebra $\cM_\cS$. In the same way the 
algebra $\cA\otimes M_N(\C)$ can be viewed as 
a $*-$subalgebra of $\cM_{\cS, N}$. More precisely, if $a\in \cA$ then the 
operator $a^\sharp$ defined by:
$$
(a^\sharp f)(x,y) := a(x+y) f(x,y), \quad \forall f\in 
B_D^2(\R^n)\otimes L^2(\R^n),
$$
belongs to $\cM$. The operator $a^\sharp$ is just the one associated with 
the zero-th order differential operator corresponding to multiplication by $a$.
The same formula allows to represent $\cA\otimes M_N(\C)$ in $\cM_{\cS,N}$.
For notional simplicity we put $N=1$ in the next result although we will use
a general $N\geq 1$ in the subsequent subsection.

\begin{prop}\
The triple $(\cA, \cM_\cS, \eth^\sharp)$ is a semifinite
 spectral triple of finite dimension equal to $n$. 
\end{prop}

\begin{proof}\
Note that the algebra $\cA$ is unital. The differential
operator $\eth$ is known to be densely defined, elliptic,
 periodic and self-adjoint on
$L^2(\R^n,\cS)$. Therefore, the operator $\eth^\sharp$ is affiliated with the
von Neumann algebra $\cm_\cS$ and it becomes self-adjoint as a densely
defined  unbounded operator on the Hilbert space $B_D^2(\R^n)\otimes
L^2(\R^n,\cS)$ with  $\eth^2 = \Delta I$ 
where $I$ is the identity operator
and $\Delta$ is the usual Laplacian.
For any smooth bounded almost periodic function $f$ on $\R^n$,
with bounded derivatives of all orders,
the commutator $[\eth, f]$ is a $0-$th order almost periodic differential 
operator and so  $[\eth^\sharp, f]$ belongs to the von Neumann algebra 
$\cm_\cS$. 

On the other hand, the pseudodifferential operator $T=(\eth^2 + I)^{-1/2}$  is 
essentially the Fourier multiplier associated with the function 
$k\mapsto \frac{1}{(\|k\|^2 +1)^{1/2}}$. Therefore, its singular 
numbers $\mu_t(T)$ can be computed explicitly as in the proof of 
Theorem \ref{Dixmier} and shown to be proportional to $t^{-1/n}$. 
\end{proof}

\subsection{Analysis of terms in the above example}

First we note that the spectral dimension is the dimension $n$
of the underlying Euclidean space and this is assumed to be odd.
It follows that the summation over $|k|$ in each term in 
the preceding theorem is
over the range $0\leq |k|\leq n-m$. 
Second we note that $m$ is always odd.

Let us write $e_1,e_2,\ldots,e_n$  for an orthonormal basis of $\R^n$,
$c(e_1),c(e_2),\ldots,c(e_n)$ for the corresponding Clifford
generators. So we have
$c(e_i)c(e_j)+c(e_j)c(e_i)=2\delta_{ij}1$
and we can write  $\eth=\sum ic(e_j)\otimes 1\partial_j$ where
1 just denotes the identity matrix.
We let $u\in \cA^\infty\otimes  \End(\C^N)$ be unitary.
Thus $u[\eth, u^*] = \sum ic(e_j)\otimes\partial_ju^* $ 
The trace is now the product of the trace on the spinor part
times the von Neumann trace composed with the matrix trace on
the matrices acting on $V$. This very simple structure enables us to
eliminate all but one of the terms in the local index formula by first
taking the trace of the product of Clifford generators. 
Note that the trace on the Clifford algebra in the spin representation
is given by
$$Tr_{Spin}(i^nc(e_1)c(e_2)\ldots c(e_n))=i^{-[(n+1)/2]}2^{(n-1)/2}$$
and the trace on any product of $0<k<n$ generators is zero.

A typical term in the local index
formula is proportional to
$$\tau_j(u[\eth,u^*]^{(k_1)}[\eth,u]^{(k_2)}
\ldots[\eth,u]^{(k_{m-1})}[\eth,u^*]^{(k_m)}
(1+\Delta)^{-|k|-m/2}) \eqno (\ast)$$
This is, up to a sign, a product of factors of the form
$(u\eth u^*-\eth)^{(k_l)}$. The Laplacian commutes with $\eth$
so that 
a typical factor  is of the form 
$\sum_i c(e_i)\otimes g_i$ and the $g_i$ are $N\otimes N$ 
matrix valued pseudodifferential operators.
Since there is always a product of an odd number of 
factors of this form $\sum_j c(e_j)\otimes g_j$ in a
term $(\ast)$ the trace on the Clifford elements will produce zero
unless $m=n$. In that case $|k|$ is forced to be zero.

Thus only one term survives in the local index theorem and that term
is (see appendix)
 $$\frac{(-1)^{(n-1)/2}}{n2^{(n-1)}}\tau_0(u[\eth,u^*][\eth,u]
[\eth,u^*][\eth,u]\ldots[\eth,u][\eth,u^*](1+\Delta)^{-n/2})$$
To compute this we
first take care of the Clifford algebra.
Using the fact that 
$[\eth,u]=-u[\eth,u^*]u$ 
we write the formula for the spectral flow as
$$\frac{(-1)^n}{n2^{(n-1)}}\tau_0(([\eth,u]u^*)^n(1+\Delta)^{-n/2})
$$

We let $Tr_N$ be the matrix trace on the auxiliary
vector space. Now $$([\eth,u]u^*)=\sum i\partial_j(u)u^*c(dx_j).
$$Writing $f_j=\partial_j(u)u^*$ we then have
$$([\D,u]u^*)^n=i^n\sum_{J=(j_1,\dots,j_n)}f_{j_1}\cdots
f_{j_n}c(e_{j_1})\cdots c(e_{j_n}),$$ where the sum is extended
over all multi-indices $J$.  Every term in the sum is a multiple
of the volume form, and so has non-zero  (spinor) trace. In terms
of permutations we have
$$([\D,u]u^*)^n=i^n\left(\sum_{\s\in\Sigma^n}(-1)^\s f_{\s(1)}\cdots
f_{\s(n)}\right) c(e_1)\cdots c(e_n)$$
$$=:\Omega i^nc(e_1)\cdots c(e_n).$$
In taking the trace we may first take the matrix trace over the
Clifford endomorphisms of the spin bundle
(with $[...]$ denoting `the integer part of')
and so, with $\tau_S=\tau\times \mbox{Tr}_N\times\mbox{Tr}_{Spin}$
\bean\tau_0\left(([\D,u]u^*)^n(1+\D^2)^{-n/2}\right)
=res_{s=0}\tau_S\left(([\D,u]u^*)^n(1+\D^2)^{-n/2-s}\right)\nno\eean 
\bean&=&res_{s=0}2^{(n-1)/2}i^{-[(n+1)/2]}(\tau\times
\mbox{Tr}_N)\left(\Omega(1+\D^2)^{-n/2-s}\right)\nno
&=&res_{s=0}
\frac{2^{(n-1)/2}}{i^{[(n+1)/2]}}
\lim\frac{1}{(2T)^n}
\int_{[-T,T]^n}\mbox{Tr}_N(\Omega)\int_{\R^n}(1+|\xi|^2)^{-n/2-s}d\xi\nno
&=&res\frac{2^{(n-1)/2}Vol(S^{n-1})}{i^{[(n+1)/2]}}
\lim\frac{1}{(2T)^n}
\int_{[-T,T]^n}\mbox{Tr}_N(\Omega)\int_0^\infty
\frac{r^{n-1}}{(1+r^2)^{n/2+s}}dr\nno
&=&res_{s=0}
\frac{2^{(n-1)/2} Vol(S^{n-1})}{i^{[(n+1)/2]}}
\lim\frac{1}{(2T)^n}
\int_{[-T,T]^n}\mbox{Tr}_N(\Omega)
\frac{\Gamma(n/2)\Gamma(s)}{2\Gamma(n/2+s)}\nno
&=&\frac{2^{(n-1)/2}}{i^{(n+1)/2}}Vol(S^{n-1})\frac{1}{2}\lim\frac{1}{(2T)^n}
\int_{[-T,T]^n}\mbox{Tr}_N(\Omega).\eean Now
$$Vol(S^{n-1})=\frac{(4\pi)^{n/2}}{2^{n-1}\Gamma(n/2)}.$$


Putting the previous calculations together gives our final result.
\begin{theorem} With the notation as above the spectral flow along 
any path joining
the Dirac operator $\eth$ to its gauge equivalent transform 
$u\eth u^*$ by an almost periodic $U(N)$ valued function
on $\R^n$ is given by the following formula:
$$sf(\eth, u\eth u^*)=  
\frac{-i^{-[(n+1)/2]}\pi^{n/2}}{\Gamma(1+n/2)2^{(n+1)/2}} 
\lim_{T\to\infty}\frac{1}{(2T)^n}\int_{(-T,T)^n}tr_N(\Omega)\ $$
\end{theorem}

\section{APPENDIX}

\subsection{Coefficients from the Local Index Theorem}

 The formula provided by the local index theorem 
for the special case considered in Section 11 is
$$ sf(\D,u^*\D u)=\frac{1}{\sqrt{2\pi i}}\sum_{m=1}^{n}
(-1)^{(m-1)/2}(\frac{(m-1)}{2})!\phi_m(u,u^*,\ldots,u, u^*)$$
where $\phi_m(u,u^*,\ldots,u, u^*)$ is given in Theorem \ref{SFLIT}.
We already know that we need only compute the top term (degree
$n$) of the local index theorem, because the Clifford trace will
kill all the other terms. Since we have a simple spectral triple,
the only multi-index $k=(k_1,...,k_n)$ which arises is zero. In
particular, we require $\alpha(0)=\frac{1}{n!}$.

The numbers $\s_{m,j}$ are defined by the equality
\ben\prod_{l=0}^{m-1}(z+l+1/2)=\sum_{j=0}^{m}z^j\s_{m,j}.\een
These are just the elementary symmetric functions of
$1/2,3/2,...,m-1/2$. When $m=0$, this is the empty product, so
$\s_{0,0}=1$.
 For $|k|=0$ we have
$h:=|k|+(n-1)/2=(n-1)/2$ and because we have simple dimension
spectrum, we only want $j=0$. Then
$\s_{(n-1)/2,0}$ is the coefficient of 1 in the product
$\prod_{l=0}^{(n-3)/2}(z+l+1/2).$ This is the product of all the
non-$z$ terms, which is
\bean(1/2)(3/2)\times\cdots\times
\left((n-3)/2+1/2\right)&=&\frac{1.
3.\cdots.(2(n-1)/2-1)}{2^{(n-1)/2}}.\eean The
reason for writing this so elaborately, is that in this form it is
obvious that it is equal to
$$\frac{1}{\sqrt{\pi}}\Gamma((n-1)/2+1/2)=\frac{1}{\sqrt{\pi}}\Gamma(n/2).$$
Combining all these calculations gives us
$$\phi_n(a_0,a_1,\dots,a_n)=\sqrt{2\pi
i}\frac{\Gamma(n/2)}{\sqrt{\pi}n!}\tau_0(a_0[\D,a_1]
\cdots[\D,a_n](1+\D^2)^{-n/2}).$$

\subsection{Constants from $Ch_n(u^*)$ and pairing} When we  pair
$\phi_n$ with the Chern character of a unitary, we divide out by
$\sqrt{2\pi i}$, which is only in the Chern character of
$(\A,\HH,\D)$ to make it compatible with the Kasparov product.
The Chern character of $u^*$ has degree $n$ component
$$ (-1)^{(n-1)/2}\left((n-1)/2)\right)!
u\otimes u^*\otimes u\otimes\cdots\otimes u^*\in\A^{\otimes
n+1}.$$ So (since $sf(\D,u\D u^*)=\frac{1}{\sqrt{2\pi
i}}\phi_n(Ch_n(u^*))$)
$$sf(\D,u\D
u^*)=\frac{(-1)^{(n-1)/2}\Gamma(n/2)\Gamma((n+1)/2)}{\sqrt{\pi}n!}
$$
$$\times \tau_0(u[\D,u^*]\cdots[\D,u^*](1+\D^2)^{-n/2})$$
Using the duplication formula for the Gamma function, we can
simplify the constant in the last displayed formula. The
duplication formula yields
$$
\Gamma(n/2)\Gamma(n/2+1/2)=\sqrt{\pi}\Gamma(n)2^{-n+1}=
\sqrt{\pi}(n-1)!2^{-n+1},$$ and inserting this gives $$ sf(\D,u\D
u^*)=\frac{(-1)^{(n-1)/2}}{n2^{(n-1)}}\tau_0(u[\D,u^*]
\cdots[\D,u^*](1+\D^2)^{-n/2}).$$

\end{document}